\documentclass[a4paper,12pt]{amsart}

\usepackage{amssymb}
\usepackage{latexsym} 
\usepackage{amsfonts} 
\usepackage{amsmath}
\usepackage{eucal} 
\usepackage{bm} 
\usepackage{bbm} 
\usepackage{graphicx} 
\usepackage[english]{varioref} 
\usepackage[nice]{nicefrac} 
\usepackage[all]{xy}
\usepackage{amsthm}

\newcommand{\supp}{\text {\rm supp}}

\newcommand{\Ad}{\text {\rm Ad}}

\def\i{^{-1}}
\def\ge{\geqslant}
\def\le{\leqslant}
\def\<{\langle} 
\def\>{\rangle}

\def\defect{\rm def}

\def\a{\alpha}

\def\g{\gamma}
\def\G{\Gamma}
\def\d{\delta}

\def\e{\epsilon}

\def\s{\sigma}
\def\t{\tau}

\def\k{\kappa}
\def\l{\lambda}

\def\Om{\Omega}

\def\bJ{\mathbb J}
\def\bS{\mathbb S}
\def\tbS{\tilde{\mathbb S}}
\def\ZZ{\mathbb Z}
\def\NN{\mathbb N}
\def\QQ{\mathbb Q}
\def\FF{\mathbb F}
\def\RR{\mathbb R}

\def\kk{\mathbf k}

\def\ca{\mathcal A}

\def\cg{\mathcal G}

\def\co{\mathcal O}
\def\cp{\mathcal P}

\def\ct{\mathcal T}
\def\cu{\mathcal U}

\def\cz{\mathcal Z}

\def\Fl{\mathcal Fl}

\def\tH{\tilde H}

\def\tx{\tilde x}
\def\ty{\tilde y}

\def\tW{\tilde W}
\def\tw{\tilde w}

\def\subset{\subseteq}

\theoremstyle{plain}
\newtheorem{thm}{Theorem}[section] 
\newtheorem*{thm*}{Theorem} 
 \newtheorem{prop}[thm]{Proposition}
 \newtheorem{lem}[thm]{Lemma}
 \newtheorem{cor}[thm]{Corollary}

\theoremstyle{definition}

\newtheorem{defn}[thm]{Definition}

\theoremstyle{remark}

\newtheorem*{rmk}{Remark}
\newtheorem*{claim*}{Claim}

\begin{document}

\author{Xuhua He}
\address{Department of Mathematics, The Hong Kong University of Science and Technology, Clear Water Bay, Kowloon, Hong Kong}
\email{maxhhe@ust.hk}
\thanks{The author is partially supported by HKRGC grant 602011.}
\title[]{Geometric and homological properties of affine Deligne-Lusztig varieties}
\keywords{affine Deligne-Lusztig varieties, $\s$-conjugacy classes}
\subjclass[2000]{14L05, 20G25}

\begin{abstract}
This paper studies affine Deligne-Lusztig varieties $X_{\tw}(b)$ in the affine flag variety of a quasi-split tamely ramified group. We describe the geometric structure of $X_{\tw}(b)$ for a minimal length element $\tw$ in the conjugacy class of an extended affine Weyl group, generalizing one of the main results in \cite{HL} to the affine case. We then provide a reduction method that relates the structure of $X_{\tw}(b)$ for arbitrary elements $\tw$ in the extended affine Weyl group to those associated with minimal length elements. Based on this reduction, we establish a connection between the dimension of affine Deligne-Lusztig varieties and the degree of the class polynomial of affine Hecke algebras. As a consequence, we prove a conjecture of G\"ortz, Haines, Kottwitz and Reuman in \cite{GHKR}. 
\end{abstract}

\maketitle

\section*{Introduction}

\subsection{} This paper discusses some geometric and homological properties of affine Deligne-Lusztig varieties in the affine flag variety of tamely ramified groups. 

To provide some context, we begin with (classical) Deligne-Lusztig varieties. Let $G$ be a connected reductive algebraic group over an algebraic closure $\kk$ of a finite field $\FF_q$, $B$ be a Borel subgroup defined over $\FF_q$, and $W$ be the associated Weyl group. We have the Bruhat decomposition $G=\sqcup_{w \in W} B \dot w B$, where $\dot w \in G$ is a representative of $w \in W$.

Let $\s$ be the Frobenius automorphism on $G$. Following \cite{DL}, the (classical) Deligne-Lusztig variety associated with $w \in W$ is a locally closed subvariety of the flag variety $G/B$ defined by $$X_w=\{g B\in G/B; g \i \s(g) \in B \dot w B\}.$$ We know that $X_w$ is always nonempty and is a smooth variety of dimension $\ell(w)$. 

The finite group $G^\s$ acts on $X_w$ and hence on its cohomology. Deligne and Lusztig showed in \cite{DL} that every irreducible representation of $G^\s$ can be realized as a direct summand of the ($l$-adic) cohomology with compact support of some Deligne-Lusztig variety, with coefficients in a local system. In \cite{Lu}, Lusztig used the cohomology to classify the irreducible representations of $G^\s$. 

\subsection{} The term ``affine Deligne-Lusztig varieties'' was first introduced by Rapoport in \cite{Ra}. Here ``affine'' refers to the fact that the notion is defined in terms of affine root systems, which arises from loop groups. 

For simplicity, we restrict our attention to the split case here. Let $G$ be a connected reductive group split over $\FF_q$ and let $L=\kk((\e))$ be the field of the Laurent series. The Frobenius automorphism $\s$ on $G$ induces an automorphism on the loop group $G(L)$, which we shall denote by the same symbol. 

Let $I$ be a $\s$-stable Iwahori subgroup of $G(L)$. By definition, the affine Deligne-Lusztig variety associated with $\tw$ in the extended affine Weyl group $\tW\cong I \backslash G(L)/I$ and $b \in G(L)$ is $$X_{\tw}(b)=\{g I \in G(L)/I; g \i b \s(g) \in I \dot \tw I\},$$ where $\dot \tw \in G(L)$ is a representative of $\tw \in \tW$. 

Understanding the emptiness/nonemptiness pattern and dimension of affine Deligne-Lusztig varieties is fundamental to understand certain aspects of Shimura varieties with Iwahori level structures. There are  two important stratifications on the special fiber of a Shimura variety: one is the Newton stratification whose strata are indexed by specific $\s$-conjugacy classes $[b] \subset G(L)$; the other is the Kottwitz-Rapoport stratification whose strata are indexed by specific elements $\tw$ of the extended affine Weyl group $\tW$. There is a close relation between the affine Deligne-Lusztig variety $X_{\tw}(b)$ and the intersection of the Newton stratum associated with $[b]$ with the Kottwitz-Rapoport stratum associated with $\tw$ (see \cite{Ha} and \cite{VW}). Our joint work with Wedhorn \cite{HW} shows that the affine Deligne-Lusztig variety is also fundamental to the study of the reduction of Shimura varieties with parahoric level structures. 

There is also an important connection between affine Deligne-Lusztig varieties and the moduli spaces of $p$-divisible groups and their analogs in function field case, the local $G$-shtukas, see \cite{HV}. 

\subsection{} The affine Deligne-Lusztig variety $X_{\tw}(b)$ for an arbitrary $\tw \in \tW$ and $b \in G(L)$ is very difficult to understand. One of the main goal of this paper is to develop a reduction method for studying the geometric and homological properties of $X_{\tw}(b)$. 

The reduction method is a combination of combinatorial, algorithmic, geometric and representation-theoretic methods. To explain why and how the method works, we first discuss the reason why the affine Deligne-Lusztig variety is more complicated than its classical (finite) counterpart. 

In the finite case, Lang's theorem implies that $G$ is a single $\s$-conjugacy class. This is the reason why a Deligne-Lusztig variety depends only on the parameter $w \in W$, with no need to choose an element $b \in G$. However, in the affine setting, the analog of Lang's theorem fails. Instead, an affine Deligne-Lusztig variety depends on two parameters: an element $\tw$ in the extended affine Weyl group and an element $b$ (or its $\s$-conjugacy class) in the loop group. Hence it is a challenging task even to describe when $X_{\tw}(b)$ is nonempty . 

\subsection{}\label{graph} To overcome this difficulty, we prove that Lang's theorem holds ``locally'' for loop groups, using a reduction method. 

We begin with the Iwahori-Bruhat decomposition $G(L)=\sqcup_{\tw \in \tW} I \dot \tw I$. Then $G(L)=\cup_{\tw \in \tW} G(L) \cdot_\s I \dot \tw I$, where $\cdot_\s$ is the $\s$-twisted conjugation action. 

Our reduction strategy can be depicted as follows. 
\[\xymatrix{\{G(L) \cdot_\s I \dot \tw I, \, \forall \tw \in \tW\} \ar@{~>}[r]^-{(1)} & \{G(L) \cdot_\s I \dot \tw I, \, \forall \tw \in \tW_{\min}\} \ar@{~>}[ld]_-{(2)}  \\ \{G(L) \cdot_\s I \dot \tw I, \forall \tw \text{ straight } \} \ar@{~>}[r]^-{(3)} & \{G(L) \cdot_\s \dot \tw, \forall \tw \text{ straight } \},}\]
where $\tW_{\min}$ is the set of elements in $\tW$ that are of minimal length in their conjugacy classes. 

In step (1), we apply  a variation of ``reduction method'' \`a la Deligne and Lusztig \cite{DL}. Our recent joint work with Nie \cite{HN} showed that minimal length elements satisfy special properties that allow us to reduce arbitrary elements to some minimal length elements. 

An element $\tw$ in $\tW$ is called straight if $\ell(\tw^n)=n \ell(\tw)$ for all $n \in \NN$, and a conjugacy class that contains a straight element is called a straight conjugacy class. It was shown in \cite{HN} that a minimal length element differs from a straight element by an element in a finite Coxeter group. Thus we may apply Lang's theorem and reduce the minimal length elements further to straight elements. This is step (2). 

Moreover, any straight element can be regarded as a basic element in the extended affine Weyl group of some Levi subgroup of $G$. Thus, using the $P$-alcove introduced in \cite{GHKR} and its generalization in \cite{GHN}, we can show that any element in $I \dot \tw I$ is $\s$-conjugate to $\dot \tw$ for some straight $\tw$. This is step (3), which completes the reduction. 

Combining this with the disjointness result in Proposition \ref{sbij}, we obtain a natural bijection between straight conjugacy classes of the extended affine Weyl group $\tW$ and $\s$-conjugacy classes of the loop group $G(L)$. This yields a new proof of Kottwitz's classification of $\s$-conjugacy classes in \cite{Ko85} and \cite{Ko97}.  

\subsection{} A similar strategy can be applied to the study of affine Deligne-Lusztig varieties, reducing the study of many geometric and homological properties of arbitrary affine Deligne-Lusztig varieties to those associated with minimal length elements. 

Although the structure of arbitrary affine Deligne-Lusztig varieties is quite complicated, the varieties associated with minimal length elements have a very nice geometric structure. We prove this in Theorem \ref{affineDL}, which generalizes one of the main theorems in \cite{HL}. 

As a consequence, our reduction method works for the homology of the affine Deligne-Lusztig varieties and we prove in Section 5 that for a superstraight element $x \in \tW$, the action of the $\s$-centralizer $\bJ_{\dot x}$ on the Borel-Moore homology of the affine Deligne-Lusztig variety $X_{\tw}(\dot x)$ factors through the action of the centralizer of $x$ in the extended affine Weyl group $\tW$. The definition of superstraight elements can be found in $\S$5.2, which includes generic elements in the cocharacter lattice of the maximal torus and superbasic element for type $A$ as special cases.  

\subsection{} The emptiness/nonemptiness pattern and dimension formula for affine Deligne-Lusztig varieties can be obtained if we can keep track of the reduction step from an arbitrary element to a minimal length element. This is accomplished via the class polynomials of affine Hecke algebras. In Section 6, we prove the ``dimension=degree'' theorem, which provides a dictionary between affine Deligne-Lusztig varieties and affine Hecke algebras in the sense that:

(1) The affine Deligne-Lusztig variety is nonempty if and only if certain class polynomial of affine Hecke algebra is nonzero. 

(2) The dimension of the affine Deligne-Lusztig variety $X_{\tw}(b)$ is equal to $\frac{1}{2} \ell(\tw)$ minus the length of the Newton point of $b$ plus a correction term given by the degree of the corresponding class polynomials.

(3) In the split case, if $b$ is superbasic, then the number of rational points in the affine Deligne-Lusztig variety $X_{\tw}(b)$ can be calculated through the corresponding class polynomial (see Section 8). 

As a consequence, we solve the emptiness/nonemptiness question for affine Deligne-Lusztig varieties in the affine Grassmannian for tamely ramified groups in Theorem \ref{Mazur}, generalizing  previous results of Rapoport-Richartz \cite{RR}, Kottwitz-Rapoport \cite{KR}, Lucarelli \cite{Luc} and Gashi \cite{Ga} for unramified groups. 

\subsection{} In the remaining sections (9--11), we study $X_{\tw}(b)$ for the case where $b$ is basic and $\tw$ is in the lowest two sided cell of $\tW$, and we prove a main conjecture of G\"ortz-Haines-Kottwitz-Reuman \cite{GHKR} and its generalization. This is achieved by combining the ``dimension=degree'' theorem with the partial conjugation method developed in \cite{He072} and the dimension formula for affine Deligne-Lusztig varieties in the affine Grassmannian in \cite{GHKR1} and \cite{Vi} 

We also give an upper bound for the dimension of $X_{\tw}(b)$ for arbitrary $\tw$ and $b$. 

\section{Group-theoretic data}

\subsection{} Let $\kk$ be an algebraic closure of a finite field $\mathbb F_q$. Let $F= \mathbb F_q( (\e))$ and $L=\kk( (\e))$ be the fields of Laurent series. Let $G$ be a connected reductive group over $F$ and splits over a tamely ramified extension of $L$. Since $\kk$ is algebraically closed, $G$ is automatically quasi-split over $L$. In this paper, we assume furthermore that the characteristic of $\kk$ does not divide the order of the fundamental group of the derived group $\pi_1(G_{der})$.\footnote{Under this assumption, the affine flag variety of $G$ is reduced (see \cite[Theorem 0.2]{PR}).}

Let $S\subset G$ be a maximal $L$-split torus defined over $F$ and $T=\cz_G(S)$ be its centralizer. Since $G$ is quasi-split over $L$, $T$ is a maximal torus. 

Let $\s$ be the Frobenius automorphism in $Gal(L/F)$, defined by $\s(\sum a_n \e^n)=\sum a_n^q \e^n$. We also denote the induced automorphism on $G(L)$ by $\s$.  

We denote by $\mathcal A$ the apartment of $G_L$ corresponding to $S$. We fix a $\sigma$-invariant alcove $\mathfrak a_C$ in $\ca$. We denote by $I \subseteq G(L)$ the Iwahori subgroup corresponding to $\mathfrak a_C$ over $L$ and by $\tbS$ the set of simple reflections at the walls of $\mathfrak a_C$. 

\subsection{}\label{splitting}
Let $N$ be the normalizer of $T$. By definition, the {\it finite Weyl group associated with $S$} is 
\[
W=N(L)/T(L)
\] 
and the {\it Iwahori-Weyl group associated with $S$} is 
\[
\tW = N(L)/T(L)_1,
\]
Here $T(L)_1$ denotes the unique parahoric subgroup of $T(L)$.

Let $\G$ be the absolute Galois group $Gal(\bar L/L)$ and $P$ be the $\G$-coinvariants of $X_*(T)$. By \cite{HR}, we may identify $T(L)/T(L)_1$ with $P$ and obtain the following important short exact sequence \[\tag{a} 0 \to P \to \tW \to W \to 0.\]

By \cite[Proposition 13]{HR}, the short exact sequence splits and we obtain a semi-direct product $\tW=P \rtimes W$ by choosing a special vertex in $\ca$. We may write an element of $\tW$ as $t^{\mu} w$ for some $\mu \in P$ and $w \in W$. 

\subsection{} Let $LG$ be the loop group associated with $G$. This is the ind-group scheme over $\kk$ which represents the functor $R \mapsto G(R((\e)))$ on the category of $\kk$-algebras. Let $\Fl=LG/{\mathbf I}$ be the {\it fppf} quotient, which is represented by an ind-scheme, ind-projective over $\kk$. We have the Iwahori-Bruhat decomposition $$G(L)=\sqcup_{\tw \in \tW} I \dot \tw I, \qquad \Fl(\kk)=\sqcup_{\tw \in \tW} I \dot \tw I/I.$$  Here $\dot \tw$ is a representative in $N(L)$ of the element $\tw \in \tW$.

Following \cite{Ra}, we define affine Deligne-Lusztig varieties as follows. 

\begin{defn}
For any $b \in G(L)$ and $\tw \in \tW$, the affine Deligne-Lusztig variety attached to $\tw$ and $b$ is the locally closed sub-ind scheme $X_{\tw}(b) \subset \Fl$ with $$X_{\tw}(b)(\kk)=\{g I \in G(L)/I; g \i b \s(g) \in I \dot \tw I\}.$$
\end{defn}

It is a finite-dimensional $\kk$-scheme, locally of finite type over $\kk$. 

For $b \in G(L)$, define the $\s$-centralizer of $b$ by $$\bJ_b=\{g \in G(L); g \i b \s(g)=b\}.$$ Then $\bJ_b$ acts on $X_{\tw}(b)$ on the left for any $\tw \in \tW$. 

\subsection{} From now on, we assume furthermore that $G$ is an adjoint group and quasi-split over $F$. As explained in \cite[Section 2]{GHN}, questions on affine Deligne-Lusztig varieties for any tamely ramified groups can be reduced to this case. 

The action of $\s$ on $G(L)$ induces an automorphism on $\tW$ and a bijection on $\tbS$. We denote the corresponding maps by $\d$. 

Moreover, we may choose a special vertex in the apartment $\ca$ such that the corresponding splitting of the short exact sequence $\S$\ref{splitting} (a) is preserved by $\d$. The induced automorphism on $W$ is also denoted by $\d$. 

Let $\Phi$ be the set of (relative) roots of $G$ over $L$ with respect to $S$ and $\Phi_a$ the set of affine roots. 
For $a \in \Phi$, we denote by $U_a \subset G$ the corresponding root subgroup and for $\a \in \Phi_a$, we denote by $U'_\a \subset G(L)$ the corresponding root subgroup scheme over $\kk$. By \cite[Section 9]{PR}, any affine real root is of the form $a+m$ for a finite root $a \in \Phi$ and $m \in \QQ$ and $U'_\a$ is one-dimensional for all affine root $\a$. Moreover, if $2 a \in \Phi$, then $m \in \ZZ$; if $\frac{1}{2} a \in \Phi$, then $m \in \frac{1}{2}+\ZZ$. We call a root $a \in \Phi$ a finite simple root if $-a$ (in $\Phi_a$) is a simple affine root with respect to $I$. Our convention here is consistent with \cite{GHKR}. 

Let $\mathbb S\subset \Phi$ be the set of simple roots. We identify $\mathbb S$ with the set of simple reflections in $W$. Then $\bS \subset \tbS$ is a $\d$-stable subset. Let $\Phi^+$ (resp. $\Phi^-$) be the set of positive (resp. negative) roots of $\Phi$. 

\subsection{} Let $G_1\subset G(L)$ be the subgroup generated by all parahoric subgroups. Set $N_1=N(L) \cap G_1$. By \cite[Proposition 5.2.12]{BT}, the quadruple $(G_1, I, N_1, \tbS)$ is a a double Tits system with affine Weyl group
\[
W_a=(N(L)\cap G_1)/(N(L)\cap I).
\] 

We may regard $W_a$ as the Iwahori-Weyl group of the simply connected cover $G_{sc}$ of the derived group $G_{der}$ of $G$. We denote by $T_{sc} \subset G_{sc}$ the maximal torus giving by the choice of $T$. Then we have the semi-direct product $$W_a=X_*(T_{sc})_\G \rtimes W.$$


By \cite{HR}, there exists a reduced root system $\Sigma$ such that \[W_a=Q^\vee(\Sigma) \rtimes W(\Sigma),\] where $Q^\vee(\Sigma)$ is the coroot lattice of $\Sigma$. In other words, we identify $Q^\vee(\Sigma)$ with $X_*(T_{\rm sc})_\Gamma$ and $W(\Sigma)$ with $W$. We simply write $Q$ for $Q^\vee(\Sigma)$. 




\subsection{} For any element $\tw \in \tW$, the length $\ell(\tw)$  is the number of ``affine root hyperplanes'' in $\ca$ separating $\tw(\mathfrak a_C)$ from $\mathfrak a_C$. 

Let $\Om$ be the subgroup of $\tW$ consisting of length $0$ elements. Then $\Om$ is the stabilizer of the base alcove $\mathfrak a_C$ in $\tW$. Let $\t \in \Om$, then for any $w, w' \in W_a$, we say that $\t w \le \t w'$ if $w \le w'$ for the Bruhat order on $W_a$. 

For any $J \subset \tbS$, let $W_J$ be the subgroup of $\tW$ generated by $s_j$ for $j \in J$ and ${}^J \tW$ (resp. $\tW^J$) be the set of minimal elements for the cosets $W_J \backslash \tW$ (resp. $\tW/W_J$). For $J, J' \subset \tbS$, we simply write ${}^J \tW^{J'}$ for ${}^J \tW \cap \tW^{J'}$. For $x \in {}^J \tW^{J'}$, we write $x(J')=J$ if $x W_{J'} x \i=W_J$. 

If $W_J$ is finite, then we denote by $w_0^J$ its longest element in $W_J$. We simply write $w_0$ for $w_0^\bS$. 

For $J \subset \bS$, let $\Phi_J$ be the set of roots spanned by $J$. Let $\Phi^+_J=\Phi_J \cap \Phi^+$ and $\Phi^-_J=\Phi_J \cap \Phi^-$. Let $M_J \subset G$ be the subgroup generated by $T$ and $U_a$ for $a \in \Phi_J$. Then $\tW_J=P \rtimes W_J$ is the Iwahori-Weyl group of $M_J$. 

\subsection{} Two elements $\tw, \tw'$ of $\tW$ are said to be $\d$-conjugate if $\tw'=\tx \tw \d(\tx) \i$ for some $\tx \in \tW$. The relation of $\d$-conjugacy is an equivalence relation and the equivalence classes are said to be $\d$-conjugacy classes. 

Let $(P/Q)_{\d}$ be the $\d$-coinvariants on $P/Q$. Let $$\kappa: \tW \to \tW/W_a \cong P/Q \to (P/Q)_{\d}$$ be the natural projection. We call $\kappa$ the {\it Kottwitz map}. 

Let $P_{\QQ}=P \otimes_\ZZ \QQ$ and $P_\QQ/W$ the quotient of $P_\QQ$ by the natural action of $W$. We may identify $P_\QQ/W$ with $P_{\QQ, +}$, where $$P_{\QQ, +}=\{\chi \in P_\QQ; \a(\chi) \ge 0, \text{ for all } \a \in \Phi^+\}.$$  Let $P_{\QQ, +}^{\d}$ be the set of $\d$-invariant points in $P_{\QQ, +}$.

Now we define a map from $\tW$ to $P_{\QQ, +}^{\d}$ as follow. For each element $\tw=t^\chi w \in \tW$, there exists $n \in \NN$ such that $\d^n=1$ and $w \d(w) \d^2(w) \cdots \d^{n-1}(w)=1$. This is because the image of $W \rtimes \<\d\>$ in $Aut(W)$ is a finite group. Then $\tw \d(\tw) \cdots \d^{n-1}(\tw)=t^\l$ for some $\l \in P$. Let $\nu_{\tw}=\l/n \in P_\QQ$ and $\bar \nu_{\tw}$ the corresponding element in $P_{\QQ, +}$. It is easy to see that $\nu_{\tw}$ is independent of the choice of $n$. 

Let $\bar \nu_{\tw}$ be the unique element in $P_{\QQ, +}$ that lies in the $W$-orbit of $\nu_{\tw}$. Since $t^\l=\tw t^{\d(\l)} \tw \i=t^{w \d(\l)}$, $\bar \nu_{\tw} \in P_{\QQ, +}^{\d}$. We call the map $\tW \to P_{\QQ, +}^{\d}, \tw \mapsto \bar \nu_{\tw}$ the {\it Newton map}. 

Define $f: \tW \to P_{\QQ, +}^{\d} \times (P/Q)_{\d}$ by $\tw \mapsto (\bar \nu_{\tw}, \kappa(\tw))$. It is the restriction to $\tW$ of the map $G(L) \to P_{\QQ, +}^{\d} \times (P/Q)_{\d}$ in \cite[4.13]{Ko97} and is constant on each $\d$-conjugacy class of $\tW$. We denote the image of the map $f$ by $B(\tW, \d)$. 

\section{Some special properties of affine Weyl groups}

In this section, we recollect some special properties on affine Weyl groups established in joint work with Nie \cite{HN}. These properties will play a crucial role in the rest of this paper. 

\subsection{}For $\tw, \tw' \in \tW$ and $i \in \tbS$, we write $\tw \xrightarrow{s_i}_\d \tw'$ if $\tw'=s_i \tw s_{\d(i)}$ and $\ell(\tw') \le \ell(\tw)$. We write $\tw \to_\d \tw'$ if there is a sequence $\tw=\tw_0, \tw_1, \cdots, \tw_n=\tw'$ of elements in $\tW$ such that for any $k$, $\tw_{k-1} \xrightarrow{s_i}_\d \tw_k$ for some $i \in \tbS$. We write $\tw \approx_\d \tw'$ if $\tw \to_\d \tw'$ and $\tw' \to_\d \tw$ and write $\tw \tilde \approx_\d \tw'$ if $\tw \approx_\d \t \tw' \d(\t) \i$ for some $\t \in \Om$. It is easy to see that $\tw \approx_\d \tw'$ if $\tw \to_\d \tw'$ and $\ell(\tw)=\ell(\tw')$. 


We call $\tw, \tw' \in \tW$ {\it elementarily strongly $\d$-conjugate} if $\ell(\tw)=\ell(\tw')$ and there exists $\tx \in \tW$ such that $\tw'=\tx \tw \d(\tx) \i$ and $\ell(\tx \tw)=\ell(\tx)+\ell(\tw)$ or $\ell(\tw \d(\tx) \i)=\ell(\tx)+\ell(\tw)$. We call $\tw, \tw'$ {\it strongly $\d$-conjugate} if there is a sequence $\tw=\tw_0, \tw_1, \cdots, \tw_n=\tw'$ such that for each $i$, $\tw_{i-1}$ is elementarily strongly $\d$-conjugate to $\tw_i$. We write $\tw \tilde \sim_\d \tw'$ if $\tw$ and $\tw'$ are strongly $\d$-conjugate. Note that $\tw \tilde \approx_\d \tw'$ implies that $\tw \tilde \sim_\d \tw'$, but the converse doesn't hold in general. 

The following result is proved in \cite[Theorem 2.10]{HN} and is a key ingredient in the reduction step (1) in $\S$\ref{graph}.

\begin{thm}\label{min} Let $\co$ be a $\d$-conjugacy class in $\tW$ and $\co_{\min}$ be the set of minimal length elements in $\co$. Then

(1) For each element $\tw \in \co$, there exists $\tw' \in \co_{\min}$ such that $\tw \rightarrow_\d\tw'$.

(2) Let $\tw, \tw' \in \co_{\min}$, then $\tw \tilde \sim_\d \tw'$.
\end{thm}

\begin{rmk}
For any $\d$-conjugacy class $\co$ of $\tW$, $\co_{\min}$ is a a single equivalence class for the relation $\tilde \sim$. However, in general $\co_{\min}$ is not a single equivalence class for the relation $\tilde \approx$. 
\end{rmk}

\subsection{} Let $\tH$ be the Hecke algebra associated with $\tW$, i.e., $\tH$ is the associated $A=\ZZ[v, v \i]$-algebra with basis $T_{\tw}$ for $\tw \in \tW$ and multiplication is given by 
\begin{gather*} T_{\tx} T_{\ty}=T_{\tx \ty}, \quad \text{ if } l(\tx)+l(\ty)=l(\tx \ty); \\ (T_s-v)(T_s+v \i)=0, \quad \text{ for } s \in \tbS. \end{gather*}

Then $T_s \i=T_s-(v-v \i)$ and $T_{\tw}$ is invertible in $H$ for all $\tw \in \tW$. 

The map $T_{\tw} \mapsto T_{\d(\tw)}$ defines an $A$-algebra automorphism of $\tH$, which we still denote by $\d$. 

Let $h, h' \in \tH$, we call $[h, h']_\d=h h'-h' \d(h)$ the {\it $\d$-commutator} of $h$ and $h'$. Let $[\tH, \tH]_\d$ be the $A$-submodule of $\tH$ generated by all $\d$-commutators. 

By Theorem \ref{min}, for any $\d$-conjugacy class $\co$ of $\tW$ and $\tw, \tw' \in \co_{\min}$, we have that $T_{\tw}\equiv T_{\tw'} \mod [\tH, \tH]_\d$. See \cite[Lemma 5.1]{HN}. 

Now for any $\d$-conjugacy class $\co$, we fix a minimal length representative $\tw_\co$. Then the image of $T_{\tw_\co}$ in $\tH/[\tH, \tH]_\d$ is independent of the choice of $\tw_\co$. It is proved in \cite[Theorem 6.7]{HN} that

\begin{thm}\label{basis}
The elements $T_{\tw_{\co}}$ form a $A$-basis of $\tH/[\tH, \tH]_\d$, here $\co$ runs over all the $\d$-conjugacy class of $\tW$. 
\end{thm}

\subsection{}\label{class} Now for any $\tw \in \tW$ and a $\d$-conjugacy class $\co$, there exists unique $f_{\tw, \co} \in A$ such that \[T_{\tw} \equiv \sum_{\co} f_{\tw, \co} T_{\tw_{\co}} \mod [\tH, \tH]_\d.\] 

By \cite[Theorem 5.3]{HN}, $f_{\tw, \co}$ is a polynomial in $\ZZ[v-v \i]$ with nonnegative coefficient and can be constructed inductively as follows. 

If $\tw$ is a minimal element in a $\d$-conjugacy class of $\tW$, then we set $$f_{\tw, \co}=\begin{cases} 1, & \text{ if } \tw \in \co \\ 0, & \text{ if } \tw \notin \co \end{cases}.$$ 

Now we assume that $\tw$ is not a minimal element in its $\d$-conjugacy class and that for any $\tw' \in \tW$ with $\ell(\tw')<\ell(\tw)$, $f_{\tw, \co}$ is constructed. By Theorem \ref{min} there exists $\tw_1 \approx_\d \tw$ and $i \in \tbS$ such that $\ell(s_i \tw_1 s_{\d(i)})<\ell(\tw_1)=\ell(\tw)$. In this case, $\ell(s_i \tw_1)<\ell(\tw)$ and we define $f_{\tw, \co}$ as $$f_{\tw, \co}=(v-v \i) f_{s_i \tw_1, \co}+f_{s_i \tw_1 s_{\d(i)}, \co}.$$

\subsection{}\label{clean} Let $V=P \otimes_\ZZ \RR$. As explained in \cite{HR}, we may regard $\Sigma$ as functions on $V$. This gives a pairing $\<, \>$ between $V$ and the root lattice of $\Sigma$. Let $\rho$ be the half sum of all positive roots in the root system $\Sigma$. 

We call an element $\tw \in \tW$ a {\it $\d$-straight element} if $\ell(\tw)=\<\bar \nu_{\tw}, 2 \rho\>$. By \cite[Lemma 1.1]{He5}, $\tw$ is $\d$-straight if and only if for any $m \in \NN$, $\ell(\tw \d(\tw) \cdots \d^{m-1}(\tw))=m \ell(\tw)$. We call a $\d$-conjugacy class {\it straight} if it contains some straight element. As we'll see later, there is a natural bijection between the set of straight $\d$-conjugacy classes of $\tW$ and the set of $\s$-conjugacy classes of loop group. 

We have the following results on straight elements and straight conjugacy classes. 

(1) The map $f: \tW \to P_{\QQ, _+}^{\d} \times (P/Q)_{\d}$ induces a bijection from the set of straight $\d$-conjugacy classes to $B(\tW, \d)$. See \cite[Theorem 3.3]{HN}. 

(2) Let $\co$ be a straight $\d$-conjugacy class of $\tW$ and $\tw, \tw' \in \co_{\min}$. Then $\tw \tilde \approx_\d \tw'$. See \cite[Theorem 3.9]{HN}. 

\

The following result \cite[Proposition 2.4 \& 2.7]{HN} relates minimal length elements with straight elements as needed in $\S$\ref{graph} step (2). 

\begin{thm}\label{min2}
Let $\co$ be a $\d$-conjugacy class of $\tW$ and $\tw \in \co$. Then there exists $\tw' \in \co_{\min}$ such that 

(1) $\tw \to_\d \tw'$;

(2) There exists $J \subset \tbS$ with $W_J$ is finite, a straight element $x \in \tW$ with $x \in {}^J \tW^{\d(J)}$ and $x \d(J)=J$ and $u \in W_J$ such that $\tw'=u x$. 
\end{thm}

\begin{rmk}
In the setting of Theorem \ref{min2} (2), we have that $f(\tw)=f(\tw')=f(x)$. See \cite[Proposition 1.2]{He5}. 
\end{rmk}

\subsection{} Let $\co$ be a $\d$-conjugacy class. Set $\nu_\co=\bar \nu_{\tw}$ for some (or equivalently, any) $\tw \in \co$ and $J_{\co}=\{i \in \bS; \<\nu_\co, \a_i\>=0\}$. The following description of straight conjugacy classes \cite[Proposition 3.2]{HN} will be used in step (3) of $\S$\ref{graph}. 

\begin{prop}\label{basicstraight}
Let $\co$ be a $\d$-conjugacy class of $\tW$. Then $\co$ is straight if and only if $\co$ contains a length $0$ element in the Iwahori-Weyl group $\tW_J=P \rtimes W_J$ for some $J \subset \bS$ with $\d(J)=J$. In this case, there exists a length $0$-element $x$ in $\tW_{J_\co}$ and $y \in W^{J_\co}$ such that $\nu_x=\nu_\co$ and $y x \d(y) \i \in \co_{\min}$. 
\end{prop}

\subsection{} Any fiber of $f: \tW \to B(\tW, \d)$ is a union of $\d$-conjugacy classes. As we'll see later, the case that the fiber is a single $\d$-conjugacy class is of particular interest. We call such $\d$-conjugacy class {\it superstraight}. A minimal length element in a superstraight $\d$-conjugacy class is called a {\it $\d$-superstraight element}. 

Let $\tw \in \Om$. Then the conjugation by $\tw$ gives a permutation on the set of simple reflections $\tbS$. We call that $\tw$ is {\it $\d$-superbasic} if each orbit of $\tw \circ \d$ on $\tbS$ is a union of connected components of the corresponding Dynkin diagram of $\tbS$. By \cite[$\S$3.5]{HN}, $\tw$ is a $\d$-superbasic element if and only if $\tW=W_1^{m_1} \times \cdots \times W_l^{m_l}$,
where $W_i$ is an affine Weyl group of type $\tilde A_{n_i-1}$ and $\tw \circ \d$ gives an order $n_i m_i$ permutation on the set of simple reflections in $W_i^{m_i}$.  

The following description of superstraight conjugacy class is obtained in \cite[Proposition 3.4]{HN}, which is analogous to Proposition \ref{basicstraight}. 

\begin{prop}\label{superstraight}
Let $\co$ be a $\d$-conjugacy class of $\tW$. Then $\co$ is superstraight if and only if there exists a $\d$-superbasic element $x$ in $\tW_{J_\co}$ and $y \in W^{J_\co}$ such that $\nu_x=\nu_\co$ and $y x \d(y) \i \in \co_{\min}$.
\end{prop}

\begin{cor}
Let $\co$ be a $\d$-straight conjugacy class of $\tW$. If $\nu_{\co}$ is regular, then $\co$ is superstraight. 
\end{cor}

\subsection{}\label{za} In particular, let $\d$ be the identity map and $\l \in P_+$ be a regular element. Let $\co=\{t^{w \l}; w \in W\}$ be the conjugacy class of $t^\l$. Then $\nu_\co=\l$ is regular and $\co$ is superstraight. Since every element in $\co$ is of the same length, $t^\l$ is a minimal length element in $\co$ and hence is a superstraight element. 

\section{$\s$-conjugacy classes}

\subsection{} We consider $\s$-conjugation action on $G(L)$, $g \cdot_\s g'=g g' \s(g) \i$. The classification of $\s$-conjugacy classes is due to Kottwitz \cite{Ko85} and \cite{Ko97}. In order to understand affine Deligne-Lusztig varieties, we also need some information on the relation between $\s$-conjugacy classes and Iwahori-Bruhat cells $I \dot \tw I$. In this section, we'll study the subsets of the form $$G(L) \cdot_\s I \dot \tw I=\{g g' \s(g) \i; g \in G(L), g' \in I \dot \tw I\}$$ for any $\tw \in \tW$. As a consequence, we obtain a new proof of Kottwitz's classification.

\begin{lem}\label{red11}
Let $\tw, \tw' \in \tW$. 

(1) If $\tw \to_\d \tw'$, then $$G(L) \cdot_\s I \dot \tw' I \subset G(L) \cdot_\s I \dot \tw I \subset G(L) \cdot_\s I \dot \tw' I \cup \cup_{\tx \in \tW, \ell(\tx)<\ell(\tw)} G(L) \cdot_\s I \dot \tx I.$$ 

(2) If $\tw \tilde \approx_\d \tw'$, then $$G(L) \cdot_\s I \dot \tw I=G(L) \cdot_\s I \dot \tw' I.$$ 
\end{lem}

Proof. It suffices to consider the case where $\tw'=\t \tw \d(\t) \i$ for some $\t \in \Om$ or $\tw \xrightarrow{i}_\d \tw'$ for some $i \in \tbS$. 

If $\tw'=\t \tw \t \i$, then $I \dot \tw' I=\dot \t I \dot \tw I \s(\dot \t) \i$ and $G(L) \cdot_\s I \dot \tw I=G(L) \cdot_\s I \dot \tw' I$. 

Now we consider the case where $\tw'=s_i \tw s_{\d(i)}$ and $\ell(\tw') \le \ell(\tw)$. By \cite[Lemma 1.6.4]{DL}, we have that $\tw=\tw'$ or $s_i \tw<\tw$ or $\tw s_{\d(i)}<\tw$. Now we prove the case where $s_{i} \tw<\tw$. The case $\tw s_{\d(i)}<\tw$ can be proved in the same way. 

Since $s_i \tw<\tw$, then $$G(L) \cdot_\s I \dot \tw I=G(L) \cdot_\s I \dot s_i I \dot s_i \dot \tw I=G(L) \cdot_\s I \dot s_i \dot \tw \s(I \dot s_i I)=G(L) \cdot_\s I \dot s_i \dot \tw I \dot s_{\d(i)} I.$$ Moreover, \[I \dot s_i \dot w I \dot s_{\d(i)} I=\begin{cases} I \dot \tw' I, & \text{ if }  \ell(\tw')=\ell(s_i \tw)+1=\ell(\tw); \\ I \dot s_i \dot \tw I \sqcup I \dot \tw' I, & \text{ if } \ell(\tw')=\ell(s_i \tw)-1=\ell(\tw)-2. \end{cases}\] In either case, $$G(L) \cdot_\s I \dot \tw' I \subset G(L) \cdot_\s I \dot \tw I \subset G(L) \cdot_\s I \dot \tw' I \cup \cup_{\tx \in \tW, \ell(\tx)<\ell(\tw)} G(L) \cdot_\s I \dot \tx I.$$ 

If moreover, $\tw \tilde \approx_\d \tw'$, then $\ell(\tw')=\ell(\tw)$. In this case, $G(L) \cdot_\s I \dot \tw I=G(L) \cdot_\s I \dot \tw' I$. \qed

\begin{lem}\label{red3} Let $J \subset \tbS$ with $W_J$ finite and $x \in {}^J \tW^{\d(J)}$ with $x \d(J)=J$. Then for any $u \in W_J$, we have that $$G(L) \cdot_\s I \dot u \dot x I=G(L) \cdot_\s I \dot x I.$$ 
\end{lem}

Proof. Let $\cp$ be the standard parahoric subgroup corresponding to $J$ and $\cu$ be its prounipotent radical. Let $\bar \cp=\cp/\cu$ be the reductive quotient of $\cp$ and $p \mapsto \bar p$ the projection from $\cp$ to $\bar \cp$. By definition, $\bar p \mapsto \dot x \overline{\s(p)} \dot x \i$ is a Frobenius morphism on $\bar \cp$. We denote it by $\s_{\dot x}$. Hence by Lang's theorem, $\bar \cp=\{ \bar p \i \dot x \overline{\s(p)} \dot x \i; p \in \cp\}$ and $$\cp \dot x/\s(\cu)=\dot x \s(\cp)/\s(\cu)=\{\bar p \i \dot x \overline{\s(p)}; p \in \cp\}.$$ Therefore, $\cp \dot x \s(\cp)=\cp \cdot_\s I \dot x I$. Similarly, $\cp \dot x \s(\cp)=\cp \cdot_\s I \dot u \dot x I$. Thus $\cp \cdot_\s I \dot u \dot x I=\cp \cdot_\s I \dot x I$ and $G(L) \cdot_\s I \dot u \dot x I=G(L) \cdot_\s I \dot x I$. \qed

\subsection{} We recall the notation of $(J, w, \d)$-alcoves in \cite[3.3]{GHN}, which generalizes the notation of $P$-alcoves introduced in \cite{GHKR}. 

\begin{defn}
Let $J \subset \bS$ with $\d(J)=J$ and $w \in W$. Let $x\in\tW$.  We say $x(\mathfrak a_C)$  is a \emph{$(J, w, \d)$-alcove}, if
\begin{enumerate}
\item[(1)] $w \i x \d(w) \in \tW_J$, and
\item[(2)] For any $a \in w(\Phi^+-\Phi^+_J)$, $\mathbf U_a \cap {}^x I \subseteq \mathbf U_a \cap I$, or equivalently, $\mathbf U_{-a} \cap {}^x I \supseteq \mathbf U_{-a} \cap I$.
\end{enumerate}
\end{defn}

Following \cite[Definition 13.1.1]{GHKR}, for $x \in \tW$, we say that $x$ is {\it $\s$-fundamental} if every element of $I \dot x I$ is $\s$-conjugate under $I$ to $\dot x$. The following result is proved in \cite[Theorem 3.3.1 \& Proposition 3.4.3]{GHN}. 

\begin{thm} \label{ghn-mainthm}
Suppose $J \subset \bS$ with $\d(J)=J$ and $w \in W$, and $x(\mathfrak a_C)$ is a $(J, w, \d)$-alcove. Set $I_M={}^{\dot w} M_J \cap I$.
Then the map
$$
\phi: I \times I_M \dot x \s(I_M) \rightarrow I \dot x I
$$
induced by $(i,m) \mapsto im \sigma(i)^{-1}$, is surjective. 

If moreover, ${}^{\dot x} \s(I_M)=I_M$, then $x$ is $\s$-fundamental.
\end{thm}

Now we describe $G(L) \cdot_\s I \dot \tw I$ for minimal length element $\tw$. 

\begin{thm}\label{smin}
If $\tw \in \tW$ is a minimal length element in its $\d$-conjugacy class, then $G(L) \cdot_\s I \dot \tw I$ is the single $\s$-conjugacy class $G(L) \cdot_\s \dot \tw$. 
\end{thm}

Proof. By Theorem \ref{min2}, $\tw \approx_\d u x$ for some $J \subset \tbS$ with $W_J$ finite,  straight element $x \in \tW$ and $u \in W_J$ with $x \in {}^J \tW^{\d(J)}$ and $x \d(J)=J$. Then by Lemma \ref{red11}, $G(L) \cdot_\s I \dot \tw I=G(L) \cdot_\s I \dot u \dot x I$. By Lemma \ref{red3}, $G(L) \cdot_\s I \dot u \dot x I=G(L) \cdot_\s I \dot x I$. 

Let $\co$ be the $\d$-conjugacy class of $x$. By Proposition \ref{basicstraight}, there exists $x_1 \in \co_{\min}$ and $y \in W^{J_{\co}}$ such that $y \i x_1 \d(y)$ is a length $0$ element in $\tW_{J_{\co}}$ and $\nu_{y \i x_1 \d(y)}=\nu_{\co}$ is dominant. Then $x_1(\mathfrak a_C)$ is a $(J_\co, y, \d)$-alcove and ${}^{\dot x_1} \s(I_M)=I_M$, where $I_M={}^{\dot y} M_{J_\co} \cap I$. Hence every element in $I \dot x_1 I$ is $\s$-conjugate under $I$ to $x_1$. 

By $\S$\ref{clean} (2), $x  \tilde \approx_\d x_1$. Now \begin{align*} G(L) \cdot_\s I \dot \tw I &=G(L) \cdot_\s I \dot x I=G(L) \cdot_\s I \dot x_1 I=G(L) \cdot_\s \dot x_1\end{align*} is a single $\s$-conjugacy class. Since $\dot \tw \in G(L) \cdot_\s I \dot \tw I$, $G(L) \cdot_\s I \dot \tw I=G(L) \cdot_\s \dot \tw$. \qed

\

We also have the following disjointness result. 

\begin{prop}\label{sbij}
Let $\tw, \tw' \in \tW$ be minimal length elements in their $\d$-conjugacy classes respectively. Then $\dot \tw$ and $\dot \tw'$ are in the same $\s$-conjugacy class of $G(L)$ if and only if $f(\tw)=f(\tw')$. 
\end{prop}

Proof. By the proof of Theorem \ref{smin}, there exists $\d$-straight elements $x, x' \in \tW$ with $f(\tw)=f(x)$, $f(\tw')=f(x')$, $\dot \tw \in G(L) \cdot_\s \dot x$ and $\dot \tw' \in G(L) \cdot_\s \dot x'$. 

If $f(\tw)=f(\tw')$, then $f(x)=f(x')$ and by $\S$\ref{clean} (1), $x$ and $x'$ are in the same $\d$-conjugacy class of $\tW$. Therefore $G(L) \cdot_\s \dot x=G(L) \cdot_\s \dot x'$. Hence $\dot \tw$ and $\dot \tw'$ are in the same $\s$-conjugacy class of $G(L)$. 

Now we prove the other direction. Notice that $\kappa(G(L) \cdot_\s \dot \tw)=\kappa(\tw)$ and $\kappa(G(L) \cdot_\s \dot \tw')=\kappa(\tw')$. Thus if $\kappa(\tw) \neq \kappa(\tw')$, then $G(L) \cdot_\s \dot \tw \cap G(L) \cdot_\s \dot \tw'=\emptyset$. 

Now assume that $\kappa(\tw)=\kappa(\tw')$ and $\bar \nu_{\tw} \neq \bar \nu_{\tw'}$. Then $\bar \nu_x \neq \bar \nu_{x'}$. If $\dot \tw$ and $\dot \tw'$ are in the same $\s$-conjugacy class of $G(L)$, then so are $\dot x$ and $\dot x'$. We may assume that $\dot x=g \dot x' \s(g) \i$ for some $g \in G(L)$, then for any $m \in \NN$, $$\dot x \s(\dot x) \cdots \s^{m-1}(\dot x)=g \dot x' \s(\dot x') \cdots \s^{m-1}(\dot x') \s^n(g) \i.$$

There exists $n \in \NN$ such that $\d^n=1$, $x \d(x) \cdots \d^{n-1}(x)=t^\mu$ and $x' \d(x') \cdots \d^{n-1}(x')=t^{\mu'}$ for some $\mu, \mu' \in P$. Then $\mu=n \nu_{x}$ and $\mu'=n \nu_{x'}$. Since $\bar \nu_{x} \neq \bar \nu_{x'}$, $\mu \notin W \cdot \mu'$. Assume that $g \in I z I$ for some $z \in \tW$. Then $I \e^{k \mu} I \cap I \dot z I \e^{k \mu'} I \dot z \i I \neq \emptyset$ for all $k \in \NN$. 

Notice that $I \dot z I \e^{k \mu'} I \dot z \i I \subset \cup_{\ty, \ty' \le z} I \ty \e^{k \mu'} (\dot \ty') \i I$. So $t^{k \mu}=\ty t^{k \mu'} (\ty') \i$ for some $\ty, \ty' \le z$. Assume that $\ty=y t^\chi$ and $\ty'=y' t^{\chi'}$ with $\chi, \chi' \in X_*(T)$ and $y, y' \in W$. Then $\ty t^{k \mu'} (\ty') \i=t^{y (k \mu'+\chi-\chi')} y (y')\i$. Hence $y=y'$ and $k \mu'+\chi-\chi'=k y \i \mu$. By definition, $\ell(t^\chi) \le \ell(\ty)+\ell(y) \le \ell(z)+\ell(w_0)$. Similarly, $\ell(t^{\chi'}) \le \ell(z)+\ell(w_0)$. Since $\mu \notin W \cdot \mu'$, then $\ell(t^{\mu'-y \i \mu}) \ge 1$. Now $$k \le \ell(t^{k(\mu'-y \i \mu)})=\ell(t^{\chi'-\chi}) \le \ell(t^{\chi'})+\ell(t^\chi) \le 2 \ell(z)+2 \ell(w_0).$$ That is a contradiction.  \qed

\

Now we obtain a new proof of Kottwitz's classification of $\s$-conjugacy classes on $G(L)$. 

\begin{thm}\label{Kotcl}
For any straight $\d$-conjugacy class $\co$ of $\tW$, we fix a minimal length representative $\tw_\co$. Then $$G(L)=\sqcup_{\co} G(L) \cdot_\s \dot \tw_\co,$$ here $\co$ runs over all the straight $\d$-conjugacy classes of $\tW$. 
\end{thm}

Proof. By Proposition \ref{sbij}, $\cup_{\co} G(L) \cdot_\s \dot \tw_\co$ is a disjoint union. Now we prove that for any $\tw \in \tW$, $I \dot \tw I \subset \sqcup_{\co} G(L) \cdot_\s \dot \tw_\co$. We argue by induction. 

If $\tw$ is a minimal length element in its $\d$-conjugacy class, then the statement follows from the proof of Theorem \ref{smin}. 

If $\tw$ is not a minimal length element in its $\d$-conjugacy class, then by Theorem \ref{min}, there exists $i \in \tbS, \tw', \tw'' \in \tW$ with $\tw' \approx_\d \tw$, $\ell(\tw'')=\ell(\tw')-2$ and $\tw' \xrightarrow{s_i}_\d \tw''$. By Lemma \ref{red11}, \begin{align*} I \dot \tw I \subset G(L) \cdot_\s I \dot \tw' I \subset G(L) \cdot_\s I \dot s_i \dot \tw' I \cup G(L) \cdot_\s I \dot \tw'' I. \end{align*}
Notice that $\ell(s_i \tw'), \ell(\tw'')<\ell(\tw')=\ell(\tw)$. Now the statement follows from  induction hypothesis on $s_i \tw'$ and $\tw''$. \qed

\subsection{} As a consequence of Theorem \ref{Kotcl}, any $\s$-conjugacy class of $G(L)$ contains a representative in $\tW$. For split groups, this was first proved in \cite[Section 7]{GHKR} in a different way.

Moreover, the map $f: \tW \to B(\tW, \d)$ is in fact the restriction of a map defined on $G(L)$, $b \mapsto (\bar \nu_b, \kappa(b))$. We call $\bar \nu_b$ the (dominant) Newton vector of $b$. 

\section{Affine Deligne-Lusztig varieties}

\subsection{} We first recall the definition and some properties on the (finite) Deligne-Lusztig varieties. 

Let $H$ be a connected reductive algebraic group over $\FF_q$ and $F$ be the Frobenius morphism on $H$. We fix an $F$-stable maximal torus $T$ and a $F$-stable Borel subgroup $B \supset T$. Let $W$ be the finite Weyl group and $\bS$ be the set of simple reflections determined by $(B, T)$. The morphism $F$ on $G$ induces bijections on $W$ and $\bS$, which we still denote by $F$. 

Following \cite[Definition 1.4]{DL}, the Deligne-Lusztig varieties associated with $w \in W$ is defined by $$X_w=X_w^H=\{g \in G/B; g \i F(g) \in B w B\}.$$ The finite group $H^F$ acts on $X_w$ in a natural way. We denote by $H^F \backslash X_w$ the space of orbits. 

It is known that $X_w$ is smooth of dimension $\ell(w)$. Moreover, we have the following result on Deligne-Lusztig varieties associated with minimal length elements, which is essentially contained in \cite[4.3]{HL}. 

\begin{thm}\label{finiteDL}
Let $w \in W$ be a minimal length element in its $F$-conjugacy class. Then $H^F \backslash X_w$ is quasi-isomorphic to the orbit space of an affine space $\kk^{\ell(w)}$ by an action of a finite torus. 
\end{thm} 

The case where $w$ is elliptic is in \cite[4.3 (a)]{HL}. The general case is not stated explicitly in \cite{HL}. However, it can be trivially reduced to the elliptic case as follows.

Let $J$ be the minimal $F$-stable subset of $\bS$ such that $w \in W_J$ and $M_J$ be the standard Levi subgroup of $H$ of type $J$. Then $w$ is an elliptic element in $W_J$ and the statement follows from the elliptic case by using the fact that $H^F \backslash X_w^H \cong M_J^F \backslash X_w^{M_J}$ (see, e.g., \cite[Section 2]{BR08}). 

\

Now we recall the ``reduction'' \`a la Deligne and Lusztig for affine Deligne-Lusztig varieties (see \cite[Proof of Theorem 1.6]{DL}, also \cite[Corollary 2.5.3]{GH}).

\begin{prop}\label{DLReduction}
Let $x\in \tW$, and let $s\in\tbS$ be a simple affine reflection.
\begin{enumerate}
\item 
If $\ell(sx\delta(s)) = \ell(x)$, then there exists a universal homeomorphism $X_x(b) \rightarrow X_{sx\delta(s)}(b)$.
\item
If $\ell(sx\delta(s)) = \ell(x)-2$, then $X_x(b)$ can be written as a disjoint union $X_x(b) = X_1 \sqcup X_2$ where $X_1$ is closed and $X_2$ is open, and such that there exist morphisms $X_1\to X_{sx\delta(s)}(b)$ and $X_2\to X_{sx}(b)$ which are compositions of a Zariski-locally trivial fiber bundle with one-dimensional fibers and a universal homeomorphism.
\end{enumerate}
\end{prop}

It is easy to see that

\begin{lem}\label{tt}
Let $\tx \in \tW$ and $\t \in \Om$. Then $X_{\tx}(b)$ is isomorphic to $X_{\t \tx \d(\t) \i}(b)$. 
\end{lem}

As a consequence of Lemma \ref{DLReduction} and Lemma \ref{tt}, we have that 

\begin{cor}\label{tttt}
Let $\tw, \tw' \in \tW$ with $\tw \tilde \approx_\d \tw'$. Then $X_{\tw}(b)$ and $X_{\tw'}(b)$ are universally homeomorphic. 
\end{cor}

\

Next we show that straight elements are automatically $\s$-fundamental. 


\begin{prop}\label{fund}
Every $\d$-straight element $x$ is $\s$-fundamental, i.e., every element of $I \dot x I$ is $\s$-conjugate under $I$ to $\dot x$. 
\end{prop}

Proof. Let $\co$ be the straight $\d$-conjugacy class of $x$. Then by the proof of Theorem \ref{smin}, there exists $x' \in \co_{\min}$ such that $x'$ is $\s$-fundamental. By $\S$\ref{clean} (2), $x \tilde \approx_\d x'$. Now the Proposition following from the following Lemmas. 

\begin{lem}
Let $x \in \tW$ and $\t \in \Om$. If $x$ is $\s$-fundamental, then so is $\t x \d(\t) \i$. 
\end{lem}

Proof. Set $x'=\t x \d(\t) \i$. Any element in $I \dot x' I=I \dot \t \dot x \s(\dot \t) \i I=\dot \t I \dot x I \s(\dot \t) \i$ is of the form $\dot \t i_1 \dot x i_2 \s(\dot \t) \i$ for some $i_1, i_2 \in I$. Since $x$ is $\s$-fundamental, there exists $i \in I$ such that $i_1 \dot x i_2=i \dot x \s(i) \i$. Thus
\begin{align*} 
\dot \t i_1 \dot x i_2 \s(\dot \t) \i &=\dot \t i \dot x \s(i) \i \s(\dot \t) \i=(\dot \t i \dot \t  \i) \dot \t \dot x \s(\dot \t) \i (\s(\dot \t) \s(i) \i \s(\dot \t) \i) \\ &=(\dot \t i \dot \t  \i) \dot \t \dot x \s(\dot \t) \i \s(\dot \t i \dot \t  \i) \i.
\end{align*}

So any element in $I \dot x' I$ is $\s$-conjugate under $I$ to $\dot \t \dot x \s(\dot \t) \i$. Hence $x'$ is $\s$-fundamental. \qed

\begin{lem}
Let $x \in \tW$ and $s \in \tbS$ with $\ell(s x \d(s))=\ell(x)$. If $x$ is $\s$-fundamental, then so is $s x \d(s)$. 
\end{lem}

Proof. We prove the case that $s x<x$. The other cases can be proved in the same way. 

Since $x$ is $\s$-fundamental, the map $I \to I \dot x I$, $g \mapsto g \dot x \s(g) \i$, is surjective. Notice that $I \dot x I \cong I \dot s I \times_I I \dot s \dot x I$. Thus the map $$I \times I \to I \dot s I \times I \dot s \dot x I, \qquad (g, g') \mapsto (g \dot s (g') \i, g' \dot s \dot x \s(g) \i)$$ is also surjective. 

The map $I \dot s I \times I \dot s \dot x I \to I \dot s \dot x I \times I \s(\dot s) I$, $(h, h') \mapsto (h', \s(h))$ is surjective. Hence the map $$I \times I \to I \dot s \dot x I \times I \s(\dot s) I, \qquad (g, g') \mapsto (g' \dot s \dot x \s(g) \i, \s(g) \s(\dot s) \s(g') \i)$$ is also surjective. 

The quotient of $I \dot s \dot x I \times I \s(\dot s) I$ by the action of the first $I$ is $I \dot s \dot x I \times_I I \s(\dot s) I \cong I \dot s \dot x \s(\dot s) I$. Hence the map $$I \mapsto I \dot s \dot x \s(\dot s) I, \qquad g' \mapsto g' \dot s \dot x \s(\dot s) \s(g') \i$$ is surjective. \qed

\

Now we prove the main result of this section, which generalizes Theorem \ref{finiteDL}. 

\begin{thm}\label{affineDL}
Let $\tw \in \tW$ be a minimal length element in its $\d$-conjugacy class. Let $b \in G(L)$ with $f(\tw)=f(b)$. Then $\dim X_{\tw}(b)=\ell(\tw)-\<\bar \nu_b, 2 \rho\>$ and $\bJ_{b} \backslash X_{\tw}(b)$ is in bijection with the orbit space of the affine space $\kk^{\ell(\tw)-\<\bar \nu_b, 2 \rho\>}$ under an action of a finite torus. 
\end{thm}

\begin{rmk}
Here we use the convention that the dimension of an empty variety is $-\infty$. 
\end{rmk}

Proof. By Theorem \ref{min2}, $\tw \approx_\d u x$ for some $J \subset \tbS$ with $W_J$ finite,  straight element $x \in \tW$ and $u \in W_J$ with $x \in {}^J \tW^{\d(J)}$ and $x \d(J)=J$. By remark of Theorem \ref{min2}, $f(\tw)=f(x)=f(b)$. Then $\dot x$ and $b$ are $\s$-conjugate and $\ell(x)=\<\bar \nu_b, 2 \rho\>$. Then $X_{\tw}(\dot x)$ and $X_{\tw}(b)$ are isomorphic and we have a natural bijection $\bJ_{\dot x} \backslash X_{\tw}(\dot x) \cong \bJ_b \backslash X_{\tw}(b)$. By Corollary \ref{tttt}, $X_{\tw}(\dot x)$ and $X_{u x}(\dot x)$ are universally homeomorphic. 

Let $\cp$ be the standard parahoric subgroup corresponding to $J$ and $\cu$ be its prounipotent radical. Let $g \in G(L)$ with $g \i \dot x \s(g) \in I \dot u \dot x I$. By the proof of Lemma \ref{red3}, there exists $p \in \cp$ such that $(g p) \i \dot x \s(g p) \in I \dot x I$. By Proposition \ref{fund}, there exists $p' \in I$ such that $(g p p') \i \dot x \s(g p p')=\dot x$. Hence $g \in \bJ_{\dot x} \cp$. 

Notice that $\bJ_{\dot x} \cp$ is in bijection with the quotient space $\bJ_{\dot x} \times_{\bJ_{\dot x} \cap \cp} \cp$. Then we have a natural bijection $$X_{u x}(\dot x) \cong \bJ_{\dot x} \times_{\bJ_{\dot x} \cap \cp} X_{u x}^{\cp}(\dot x),$$ where $X_{u x}^{\cp}(\dot x)=\{g \in \cp/I; g \i \dot x \s(g) \in I \dot u \dot x I\}$. In particular, $\bJ_b \backslash X_{\tw}(b) \cong (\bJ_{\dot x} \cap \cp) \backslash X_{u x}^{\cp}(\dot x)$. 

Moreover, the projection $\Fl \to LG/\cp$ sends $X_{u x}(\dot x)$ onto $\bJ_{\dot x}/(\bJ_{\dot x} \cap \cp)$ and each fiber is isomorphic to $X_{u x}^{\cp}(\dot x)$. Hence $\dim X_{\tw}(b)=\dim X_{u x}(\dot x)=\dim X_{u x}^{\cp}(\dot x)$. 

Let $\s_{\dot x}$ be the Frobenius morphism on $\bar \cp$ defined by $\bar p \mapsto \dot x \overline{\s(p)} \dot x \i$. Notice that $\bar I=I/\cu$ is a Borel subgroup of $\bar \cp$. Therefore $X_{u x}^{\cp}(\dot x)$ is isomorphic to the classical Deligne-Lusztig variety $$X'=\{\bar p \in \bar \cp/\bar I; \bar p \i \s_{\dot x}(\bar p) \in \bar I \dot u \bar I\}/\bar I$$ of $\bar \cp$. The action of $\bJ_{\dot x} \cap \cp$ on $X_{u x}^{\cp}(\dot x)$ factors through an action of $(\bJ_{\dot x} \cap \cp) \cu/\cu \cong \bar \cp^{\s_{\dot x}}$. Hence we have a natural bijection $(\bJ_{\dot x} \cap \cp) \backslash X_{u x}^{\cp}(\dot x) \cong \bar \cp^{\s_{\dot x}} \backslash X'$.

The map $v \mapsto v x$ sends a $\s_{\dot x}$-conjugacy class in $W_J$ into a $\d$-conjugacy class in $\tW$. Since $u x$ is a minimal length element in its $\d$-conjugacy class, $u$ is a minimal length element in its $\s_{\dot x}$-conjugacy class in $W_J$. Moreover, $\ell(u)=\ell(u x)-\ell(x)=\ell(\tw)-\<\bar \nu_b, 2 \rho\>$. Now the Theorem follows from Theorem \ref{finiteDL}. \qed

\

By the same argument, we have the following result. 

\begin{prop}\label{fred}
Let $J \subset \tbS$ with $W_J$ finite, $x \in {}^J \tW^{\d(J)}$ with $x \d(J)=J$ and $u \in W_J$. Then for any $b \in G(L)$, $$\dim X_{u x}(b)=\dim X_x(b)+\ell(u).$$
\end{prop}\

\section{homology of affine Deligne-Lusztig varieties}

\subsection{} Notice that $X_{\tw}(b)=\lim\limits_{\rightarrow} X_i$ for some closed subschemes $X_1 \subset X_2 \subset \cdots \subset X_{\tw}(b)$ of finite type. Let $l$ be a prime not equal to the characteristic of $\kk$. Then $H^j_c(X_i, \bar \QQ_l)$ is defined for all $j$. Set $$H^{BM}_j(X_{\tw}(b), \bar \QQ_l)=\lim\limits_{\rightarrow} H^j_c(X_i, \bar \QQ_l)^*.$$ Then $H^{BM}_j(X_{\tw}(b), \bar \QQ_l)$ is a smooth representation of $\bJ_b$. Hence it is a semisimple module for any open compact subgroup of $\bJ_b$. 

The following result can be proved along the line of \cite[Theorem 1.6]{DL}. 

\begin{lem}\label{inddd}
Let $b \in G(L)$ and $K$ be an open compact subgroup of $\bJ_b$. Let $\tw \in \tW$, and let $i \in\tbS$. Then 

(1) If $\ell(s_i \tw s_{\d(i)})=\ell(\tw)$, then for any $j \in \ZZ$, $$H^{BM}_j(X_{\tw}(b), \bar \QQ_l) \cong H^{BM}_j(X_{s_i \tw s_{\d(i)}}(b), \bar \QQ_l)$$ as $\bJ_b$-modules. 

(2) If $\ell(s_i \tw s_{\d(i)}) =\ell(\tw)-2$, then for any simple $K$-module $M$ that is a direct summand of $\oplus_j H^{BM}_j(X_{\tw}(b), \bar \QQ_l)$, $M$ is also a direct summand of $\oplus_j H^{BM}_j(X_{s_i \tw s_{\d(i)}}(b), \bar \QQ_l) \oplus \oplus_j H^{BM}_j(X_{s_i \tw}(b), \bar \QQ_l)$. 
\end{lem}

\

Now we can show the following ``finiteness'' result. 

\begin{thm}\label{BM}
Let $b \in G(L)$ and $K$ be an open compact subgroup of $\bJ_b$. Let $M$ be a simple $K$-module. If $M$ is a direct summand of $\oplus_{\tw \in \tW} \oplus_j H^{BM}_j(X_{\tw}(b), \bar \QQ_l)$, then $M$ is a direct summand of $H^{BM}_j(X_{\tx}(b), \bar \QQ_l)$ for some $j \in \ZZ$ and $\tx$ is a minimal length element in its $\d$-conjugacy class and $f(\tx)=f(b)$. 
\end{thm}

\begin{rmk}
Given $b$, there are only finitely many $\tx$ satisfying the conditions above. 
\end{rmk}

Proof. Let $\tw \in \tW$ be a minimal length element in its $\d$-conjugacy class such that $M$ is a direct summand of $\oplus_j H^{BM}_j(X_{\tw}(b), \bar \QQ_l)$. By Theorem \ref{min} and Lemma \ref{inddd}, $\tw$ is a minimal length element in its $\d$-conjugacy class. By Theorem \ref{smin}, $G(L) \cdot_\s I \dot \tw I$ is a single $\s$-conjugacy class. Hence $X_{\tw}(b) \neq \emptyset$ if and only if $b \in G(L) \cdot_\s I \dot \tw I$, i.e. $f(b)=f(\tw)$. \qed 

\

In the rest of this section, we discuss the special case where $b$ is a superstraight element. We first describe its centralizer in the affine Weyl group and in the loop group. 

\begin{prop}\label{Jx}
Let $J \subset \bS$ and $x$ be a $\d$-superbasic element in $\tW_J$ with $\nu_x \in P_{\QQ, +}$ and $J=\{i \in \bS; \<\nu_x, \a_i\>=0\}$. 

(1) Let $Z_{\tW, \d}(x)=\{w \in \tW; w \i x \d(w)=x\}$ be the $\d$-centralizer of $x$ in $\tW$. Then $Z_{\tW, \d}(x)$ consists of length $0$ elements $y$ in $\tW_J$ with $y \i x \d(y)=x$. 

(2) $M_J(L) \cap I \cap \bJ_{\dot x}$ is a normal subgroup of $\bJ_{\dot x}$ and $$\bJ_{\dot x}/\bigl(M_J(L) \cap I \cap \bJ_{\dot x} \bigr) \cong Z_{\tW, \d}(x).$$ 
\end{prop}

Proof. Let $n \in \NN$ with $\d^n=id$ and $x \d(x) \cdots \d^{n-1}(x)=t^{n \nu_x}$. Hence $\dot x \s(\dot x) \cdots \s^{n-1}(\dot x) \in \e^{n \nu_x} T(L)_1$. Thus after replacing $\dot x$ by $h \i \dot x \s(h)$ for a suitable $h \in T(L)_1$, we may assume that $\dot x \s(\dot x) \cdots \s^{n-1}(\dot x)=\e^{n \nu_x}$. 

(1) Let $w \in \tW$ with $w \i x \d(w)=x$. Then \begin{align*} w \i t^{n \nu_x} w &=w \i x \d(x) \cdots \d^{n-1}(x) \d^n(w) \\ &=(w \i x \d(w)) \d(w \i x \d(w)) \cdots \d^{n-1}(w \i x \d(w)) \\&=x \d(x) \cdots \d^{n-1}(x)=t^{n \nu_x}. \end{align*}

Thus $w \nu_x=\nu_x$ and $w \in \tW_J$. 

Let $\tilde \bS'$ be the set of simple reflection in $\tW_J$. Let $\ell_J$ be the length function and $<_J$ be the Bruhat order on $\tW_J$. If $\ell_J(w)>0$, then there exists $i \in \tilde \bS'$ such that $s_i w<_J w$. Since the map $\tW_J \to \tW_J$, $y \mapsto x \d(y) x \i$ preserve the Bruhat order, we have that $$x s_{\d(i)} x \i w=x s_{\d(i)} \d(w) x \i<_J x \d(w) x \i=w.$$ Similarly, $(\Ad(x) \circ \d)^m(s_i) w<_J w$ for all $m \in \ZZ$. However, $x$ is $\d$-superbasic for $\tW_J$ and thus $\tilde \bS'$ is a single orbit of $\Ad(x) \circ \d$. Hence $s_j w<w$ for all $j \in \tilde \bS'$. This is impossible. Thus $\ell_J(w)=0$.

(2) Let $g \in \bJ_{\dot x}$. Then $g \i \dot x \s(g)=\dot x$. We prove that 

(a) $g \in M_J(L)$. 

Similar to part (1), we have that $g \i \e^{n \nu_x} \s^n(g)=\e^{n \nu_x}$ and $\e^{-n \nu_x} g \e^{n \nu_x}=\s^n(g)$. Let $U(L)$ be the subgroup generated by $U_a$ for $a \in \Phi^+-\Phi^+_J$. We may assume that $g=u m \dot x m' u'$ for $u, u' \in U(L)$, $m, m' \in M_J(L)$ and $w \in {}^J W^J$. Since $\nu_x$ is neutral for $M_J(L)$, we have that $$(\e^{-n \nu_x} u \e^{n \nu_x}) m (\e^{-n \nu_x} \dot w \e^{n \nu_x}) m' (\e^{-n \nu_x} u' \e^{n \nu_x})=\s^n(u) \s^n(m) \s^n(\dot w) \s^n(m') \s^n(u').$$

Therefore $m (\e^{-n \nu_x} \dot w \e^{n \nu_x}) m'=\s^n(m \dot w m')$. 

Let $M'(L)$ be the derived group of $M_J(L)$ and $Z$ be the center of $M_J(L)$. Then $Z \subset T(L)$ and $M_J(L)=M'(L) Z$. We may write $m$ as $m_1 z$ and $m'$ as $m'_1 z'$ for $m, m' \in M'$ and $z, z' \in Z$. Then we have that $$z (\e^{-n \nu_x} \dot w \e^{n \nu_x}) z'=z \e^{n (w \nu_x-\nu_x)} (\dot w z' \dot w \i) \dot w=\s^n(z \dot w z').$$ 

Notice that $z \dot w z' \dot w \i \in T(L)$. Hence $$\e^{n(w \nu_x-\nu_x)}=(z \dot w z' \dot w \i) \i \s^n(z \dot w z' \dot w \i) (\s^n(\dot w) \dot w \i).$$ We have that $(z \dot w z' \dot w \i) \i \s^n(z \dot w z' \dot w \i), \s^n(\dot w) \dot w \i \in T(L)_1$. Thus $\e^{n(w \nu_x-\nu_x)} \in T(L)_1$ and $w \nu_x=\nu_x$. Therefore $\<\nu_x, w \i \a\>=\<w \nu_x, \a\>=\<\nu_x, \a\>=0$ for all $\a \in \Phi_J$. Therefore $w \i \a \in \Phi_J$ for all $\a \in \Phi_J$ and $w=1$. Hence $g=u m$ for some $u \in U(L)$ and $m \in M_J(L)$ and $\e^{-n \nu_x} u \e^{n \nu_x}=\s^n(u)$. Notice that $\<\nu_x, \a\>>0$ for all $\a \in \Phi^+-\Phi^+_J$, we must have that $u=1$. 

Thus $g \in M_J(L)$ and (a) is proved. 

Now set $I'=I \cap M_J(L)$. Then we may assume that $g \in I' \dot y I'$ for some $y \in \tW_J$. Since $x$ is a length $0$ element in $\tW_J$, $\dot x I'=I' \dot x$. Thus $g=\dot x \s(g) \dot x \i \in \dot x I' \s(\dot y) I'  \dot x \i=I' \dot x \s(\dot y) \dot x \i I'$ and $I' \dot y I' \cap I' \dot x \s(\dot y) \dot x \i I' \neq \emptyset$. Therefore $y=x \d(y) x \i$ and $y \in Z_{\tW, \d}(x)$ is a length $0$ element in $\tW_J$. Thus $g \in I' \dot y$ and $$g \i (\bJ_{\dot x} \cap I') g \subset \bJ_{\dot x} \cap g \i I' g=\bJ_{\dot x} \cap \dot y \i I' \dot y=\bJ_{\dot x} \cap I'.$$

Hence $\bJ_{\dot x} \cap I'$ is a normal subgroup of $\bJ_{\dot x}$. Moreover, we have an injective group homomorphism $$\pi: \bJ_{\dot x}/(\bJ_{\dot x} \cap I') \to Z_{\tW, \d}(x).$$
On the other hand, let $y \in Z_{\tW, \d}(x)$. Then $y \i x \d(y)=x$ and $\dot y \i \dot x \s(\dot y) \in \dot x T(L)_1$. Replacing $\dot y$ by $\dot y h$ for a suitable element $h \in T(L)_1$, we have that $\dot y \i \dot x \s(\dot y)=\dot x$. Hence the map $\pi$ is surjective and $\bJ_{\dot x}/(\bJ_{\dot x} \cap I') \cong Z_{\tW, \d}(x)$. \qed

\begin{thm}\label{cct}
Let $x \in \tW$ be a $\d$-superstraight element. Then for any $\tw \in \tW$ with $X_{\tw}(\dot x) \neq \emptyset$ and $j \in \ZZ$, the action of $\bJ_{\dot x}$ on $H^{BM}_j(X_{\tw}(b), \bar \QQ_l)$ factors through an action of the $\d$-centralizer $Z_{\tW, \d}(x)$ of $x$. 
\end{thm}

Proof.  By Proposition \ref{superstraight}, there exists a superbasic element $x_1 \in \tW_{J_\co}$ and $y \in W^{J_\co}$ such that $\nu_{x_1}=\nu_\co$ and $y x_1 \d(y) \i \in \co_{\min}$. 

After $\s$-conjugating by a suitable element in $G(L)$, we may assume that $\dot x=\dot y \dot x_1 \s(\dot y) \i$. Then $\bJ_{\dot x}=\dot y \bJ_{\dot x_1} \dot y \i$. 

By Proposition \ref{Jx}, $M_{J_\co}(L) \cap I \cap \bJ_{\dot x_1}$ is a normal subgroup of $\bJ_{\dot x_1}$ and \[\tag{a} \bJ_{\dot x_1}/(M_{J_\co}(L) \cap I \cap \bJ_{\dot x_1}) \cong Z_{\tW, \d}(x_1).\] Set $K=\dot y(M_{J_\co}(L) \cap I) \dot y \i \cap \bJ_{\dot x}$. Then $K$ is a normal subgroup of $\bJ_{\dot x}=\dot y \bJ_{\dot x_1} \dot y \i$ and $\bJ_{\dot x}/K \cong y Z_{\tW, \d}(x_1) y \i=Z_{\tW, \d}(x)$. 

By the proof of Theorem \ref{affineDL}, \[\tag{b} X_{x}(\dot x)=\bJ_{\dot x} I/I \cong \bJ_{\dot x}/(\bJ_{\dot x} \cap I).\] Since $y \in W^{J_\co}$, $\dot y (M_{J_\co}(L) \cap I) \dot y \i \subset I$. Thus $K$ acts trivially on $X_{x}(\dot x)$ and also trivially on $H^{BM}_j(X_{x}(\dot x), \bar \QQ_l)$ for any $j \in \ZZ$. By Theorem \ref{BM} and Lemma \ref{inddd} (1), any simple $K$-module that appears as a direct summand of $H^{BM}_j(X_{\tw}(\dot x), \bar \QQ_l)$ is trivial. Hence $K$ acts trivially on $H^{BM}_j(X_{\tw}(\dot x), \bar \QQ_l)$ and the action of $\bJ_{\dot x}$ factor through $Z_{\tw, \d}(x)$.  \qed 

\subsection{}\label{5.2} Now we discuss some special cases. 

Let $\l \in P^\d$ be a regular element and $\co$ be the $\d$-conjugacy class that contains $t^\l$. Then $\nu_\co=\l$ is regular. By $\S$\ref{za},  $\co$ is superstraight. By $\S$\ref{clean}, $t^\l$ is $\d$-superstraight since $\ell(t^\l)=\<\l, 2 \rho\>$. By Proposition \ref{Jx}, $Z_{\tW, \d}(t^\l)=\{t^\mu; \mu=\d(\mu) \in P\} \cong P^\d$. Therefore by Theorem \ref{cct}, for any $\tw \in \tW$ and $j \in \ZZ$, the action of $\bJ_{\e^\l}$ on $H^{BM}_j(X_{\tw}(\e^\l), \bar \QQ_l)$ factors through an action of $P^\d$. The special case for split $SL_2$ and $SL_3$ was first obtained by Zbarsky \cite{Zb} via direct calculation. 

Let $\t$ be a $\d$-superbasic element. Then $\t$ is $\d$-superstraight and by Proposition \ref{Jx}, $Z_{\tW, \d}(\t)=\Om^\d$. Therefore by Theorem \ref{cct}, for any $\tw \in \tW$ and $j \in \ZZ$, the action of $\bJ_{\dot \t}$ on $H^{BM}_j(X_{\tw}(\dot \t), \bar \QQ_l)$ factors through an action of $\Om^\d$. 

\section{Dimension ``='' degree} 

In this section, we give a formula which relates the dimension of affine Delgine-Lusztig varieties with the degree of the class polynomials. 

\begin{thm}\label{class3}
Let $b \in G(L)$ and $\tw \in \tW$. Then 

\[\dim (X_{\tw}(b))=\max_{\co} \frac{1}{2}(\ell(\tw)+\ell(\co)+\deg(f_{\tw, \co}))-\<\bar \nu_b, 2\rho\>,\] here $\co$ runs over $\d$-conjugacy class of $\tW$ with $f(\co)=f(b)$ and $\ell(\co)$ is the length of any minimal length element in $\co$.
\end{thm}

\begin{rmk}
Here we use the convention that the dimension of an empty variety and the degree of a zero polynomial are both $-\infty$. 
\end{rmk}

Proof. We argue by induction on $\ell(\tw)$. 

Let $\co'$ be the $\d$-conjugacy class that contains $\tw$. 

If $\tw \in \co'_{\min}$, then 
\[\frac{1}{2}(\ell(\tw)+\ell(\co)+\deg(f_{\tw, \co}))=\begin{cases} \ell(\tw), & \text{ if } \co=\co'; \\ -\infty, & \text{ if } \co \neq \co'. \end{cases}\]

Thus the right hand side of $(a)$ is $-\infty$ if $f(\tw) \neq f(b)$ and $\ell(\tw)-\<\bar \nu_b, 2 \rho\>$ if $f(\tw)=f(b)$. Now the statement follows from Theorem \ref{affineDL}. 

If $\tw \notin \co'_{\min}$, by Theorem \ref{min} there exists $\tw_1 \approx_\d \tw$ and $i \in \tbS$ such that $\ell(s_i \tw_1 s_{\d(i)})<\ell(\tw_1)=\ell(\tw)$. Then by $\S$\ref{class}, $$f_{\tw, \co}=(v-v \i) f_{s_i \tw_1, \co}+f_{s_i \tw_1 s_{\d(i)}, \co}.$$ 

Notice that the coefficients of class polynomials are nonnegative integers. Thus the highest term of $(v-v \i) f_{s_i \tw_1, \co}$ and $f_{s_i \tw_1 s_{\d(i)}, \co}$ do not cancel with each other. Hence $\deg(f_{\tw, \co})=\max\{\deg(f_{s_i \tw_1, \co})+1, \deg(f_{s_i \tw_1 s_{\d(i)}, \co})\}$ and $$\tiny \ell(\tw)+\deg(f_{\tw, \co})=\max\{\ell(s_i \tw_1)+\deg(f_{s_i \tw_1, \co}), \ell(s_i \tw_1 s_{\d(i)})+\deg(f_{s_i \tw_1 s_{\d(i)}, \co})\}+2.$$ By Lemma \ref{DLReduction}, $$\dim(X_{\tw}(b))=\dim(X_{\tw_1}(b))=\max\{\dim(X_{s_i \tw_1, \d_F}(b)), \dim(X_{s_i \tw_1 s_{\d(i)}}(b))\}+1.$$ Now the statement follows from induction hypothesis. 
\qed 

\subsection{} This is a key part of this paper. On one hand, it provides both theoretic and practical way to determine the dimension of affine Deligne-Lusztig varieties. On the other hand, it shows that the dimension and nonemptiness pattern of affine Deligne-Lusztig varieties $X_{\tw}(b)$ only depend on the data $(\tW, \d, \tw, f(b))$ and thus independent of the choice of $G$. 

Now we construct another group (split over $L$) with the same data $(\tW, \d)$ as follows. 

Recall that $G$ is an adjoint group. Let $H$ be a connected semisimple group of adjoint type over $\kk$ with the root system $\Sigma$. The Iwahori-Weyl group of $H(L)$ is just $\tW$. The standard Frobenius automorphism on $H$ induces an automorphism on $H(L)$, which we denote by $\s_0$. The corresponding automorphism on $\tW$ is just the identity map. 

If $\d$ is not identity, then $\d$ induces a nontrivial diagram automorphism on the (finite) Dynkin diagram of $H$. Since $H$ is adjoint, such automorphism can be lifted to an automorphism (of algebraic groups) on $H$. We denote it by $\t$. The automorphism $\t \circ \s_0$ on $H(L)$ induces the desired automorphism $\d$ on $\tW$. 

The question on the dimension and nonemptiness of affine Deligne-Lusztig varieties for $G$ can be reduced to the same question for the group $H(L)^{\t \circ \s_0}$. 

In other words, for dimension/nonemptiness pattern of affine Deligne-Lusztig varieties, it suffices to consider the case where $G$ is split over $L$. 

\

Now we list some easy consequences of Theorem \ref{class3}. 

\begin{cor}\label{6.2}
Let $b \in G(L)$ and $\tw \in \tW$. Then $X_{\tw}(b) \neq \emptyset$ if and only if $f_{\tw, \co} \neq 0$ for some $\d$-conjugacy class $\co$ of $\tW$ with $f(\co)=f(b)$.
\end{cor}

\begin{cor}
Let $\co$ be a superstraight $\d$-conjugacy class in $\tW$ and $x \in \co_{\min}$ is a $\d$-superstraight element. Then for any $\tw \in \tW$, \[\dim (X_{\tw}(\dot x))=\frac{1}{2}(\ell(\tw)+\deg(f_{\tw, \co}))-\<\nu_{\co}, \rho\>.\]

In particular, $X_{\tw}(\dot x) \neq \emptyset$ if and only if $f_{\tw, \co} \neq 0$. 
\end{cor}

\subsection{} By Theorem \ref{class3} and its consequences, we obtain all the information on the emptiness/nonemptinss and dimension formula of affine Deligne-Lusztig varieties if the degree of the class polynomials are known. The latter problem requires a thorough understanding of trace formula of all finite dimensional representation of affine Hecke algebras and thus computing class polynomials is quite hard in general and not known yet. 

However, at present, Theorem \ref{class3} is still quite useful in the study of affine Deligne-Lusztig varieties. It implies that the emptiness/nonemptinss and dimension formula of affine Deligne-Lusztig varieties only rely on the reduction method. Such observation will play a key role in the proof of emptiness/nonemptiness pattern of affine Deligne-Lusztig varieties for basic $b$, which will be discussed in \cite{GHN}. 

\subsection{} One can show that the reduction method also works in the $p$-adic case and the $p$-adic variant of $X_{\tw}(b)$ (for $b \in \tW$) is nonempty if and only if $f_{\tw, \co} \neq 0$ for some $\d$-conjugacy class $\co$ of $\tW$ with $f(\co)=f(b)$. 

There is no known good notion of dimension for the $p$-adic variant of affine Deligne-Lusztig varieties. However, one may hope that once it is established, then the dimension of $X_{\tw}(b)$ should agree in the $p$-adic and function field case, and thus the ``dimension=degree'' theorem remains valid for $p$-adic case. 

\section{Mazur's inequality and its converse}

\subsection{} In this section, we discuss some application of Theorem \ref{class3} to the hyperspecial case. 

Let $\cg$ be the smooth affine group scheme associated with the special vertex of the Bruhat-Tits building of $G$ and $L^+\cg$ be the infinite-dimensional affine group scheme defined by $L^+\cg(R)=\cg(R[[\e]])$. The fpqc quotient $Gr=LG/L^+\cg$ is called the {\it (twisted) affine Grassmannian}. We have the Cartan decomposition (see \cite{Ri}). 
\begin{gather*} G(L)=\sqcup_{\mu \in P_+} L^+\cg(\kk) \e^\mu L^+\cg(\kk), \\ Gr(\kk)=\sqcup_{\mu \in P_+} L^+\cg(\kk) \e^\mu L^+\cg(\kk)/L^+\cg(\kk).\end{gather*}

The affine Deligne-Lusztig variety in the affine Grassmannian associated with $b \in G(L)$ and $\mu \in P_+$ is the locally closed sub-ind scheme $X_\mu(b) \subset Gr$ with $$X_\mu(b)(\kk)=\{g \in G(L); g \i b \s(g) \in L^+\cg(\kk) \e^\mu L^+\cg(\kk)\}/L^+\cg(\kk) \subset Gr(\kk).\footnote{Here we distinguish between the element $\mu \in P$ and the element $t^\mu \in \tW$. In particular, $X_\mu(b)$ is a sub-ind scheme of the affine Grassmannian and $X_{t^\mu}(b)$ is a sub-ind scheme of the affine flag variety.}$$

\subsection{} Let $J \subset \bS$ with $\d(J)=J$. Let $X_J$ be the quotient of $P$ by the sublattice of $Q$ generated by the simple coroots in $J$. The action of $\d$ extends in a natural way to $X_J$. Let $Y_J=X_J/(1-\d)X_J$ be the coinvariants of this action and $Y_J^+$ be the image of $P_+$ in $Y_J$. 

Following \cite{Ga}, we define a partial order $\preceq_J$ on $Y_J^+$ as follows. For $\mu, \mu' \in Y_J$, we write $\mu \preceq_J \mu'$ if $\mu'-\mu$ is a nonnegative integral linear combination of the image in $Y_J$ of the simple coroots in $\bS-J$. 

Now we have the following result. 

\begin{thm}\label{Mazur}
Let $J \subset \bS$ with $\d(J)=J$. Let $\mu \in P_+$ and $b \in M_J(L)$ be a basic element such that $\kappa_{M_J}(b) \in Y^+_J$. Then $X_\mu(b) \neq \emptyset$ if and only if $\kappa_{M_J}(b) \preceq_J \mu$. 
\end{thm}

\begin{rmk}
The ``only if'' side was proved by Rapoport and Richartz in \cite{RR}, which generalize Mazur's theorem . The ``if'' side was conjectured by Kottwitz and Rapoport and proved for type A and C in \cite{KR}. It was then proved by Lucarelli \cite{Luc} for classical split groups and then by Gashi \cite{Ga} for unramified cases.
\end{rmk}

Proof. Since $L^+\cg(\kk) t^\mu L^+\cg(\kk)=\sqcup_{\tw \in W t^\mu W} I \dot \tw I$, $X_\mu(b) \neq \emptyset$ if and only if $X_{\tw}(b) \neq \emptyset$ for some $\tw \in W t^\mu W$. Let $C$ be the set of $\d$-conjugacy classes $\co$ in $\tW$ with $f(\co)=f(b)$. Then $C$ only depends on $\kappa_{M_J}(b)$. Thus by Corollary \ref{6.2}, 

(a) $X_\mu(b) \neq \emptyset$ if and only if $f_{\tw, \co} \neq 0$ for some $\tw \in W t^\mu W$ and $\co \in C$. 

The latter one only depends on the combinational data $(\tW, \d, \mu, J, \kappa_{M_J}(b))$ and is independent of the loop group $G$. 

To the pair $(\tW, \d)$, we may associate a quasi-split unramified group $H$. By \cite[Theorem 5.1]{Ga}, the statement holds for $H$. Thus by (a), we have that 

(b) $\kappa_{M_J}(b) \preceq_J \mu$ if and only if $f_{\tw, \co} \neq 0$ for some $\tw \in W t^\mu W$ and $\co \in C$. 

The theorem then follows from (a) and (b). \qed

\section{Superbasic elements}

For $b \in G(L)$, if $X_{\tw}(b)$ is nonempty, then usually it has infinitely many irreducible components. However, for superbasic elements, the affine Deligne-Lusztig variety is much nicer. 

\begin{prop}
Assume that $G$ is semisimple and $x \in \tW$ is a $\d$-superbasic element. Then for any $\tw \in \tW$, $X_{\tw}(\dot x)$ has only finitely many irreducible components. 
\end{prop}

Proof. We argue by induction on $\ell(\tw)$. 

Let $\co$ be the $\d$-conjugacy class containing $\tw$. 

If $\tw \in \co_{\min}$, then by Theorem \ref{smin} and Proposition \ref{sbij}, $X_{\tw}(\dot x) \neq \emptyset$ if and only if $f(\tw)=f(x)$. Since $x$ is $\d$-superbasic, we have that $\tw$ is $\d$-conjugate to $x$. By the proof of Theorem \ref{cct}, $X_{\tw}(\dot x) \cong X_{\tw}(\dot \tw)$ can be identified with a subset of $\Om$. Hence $X_{\tw}(\dot x)$ is a finite set. 

If $\tw \notin \co_{\min}$, then by Theorem \ref{min} there exists $\tw_1 \approx_\d \tw$ and $i \in \tbS$ such that $\ell(s_i \tw_1 s_{\d(i)})<\ell(\tw_1)=\ell(\tw)$. By induction hypothesis $X_{s_i \tw_1 s_{\d(i)}}(\dot x)$ and $X_{s_i \tw_1}(\dot x)$ both have finitely many irreducible components. Therefore by Proposition \ref{DLReduction}, $X_{\tw}(\dot x)$ is universal homeomorphic to $X_{\tw_1}(\dot x)$ and has only finitely many irreducible components. 
\qed

\begin{cor}
Assume that $G$ is semisimple and $x \in \tW$ is a $\d$-superbasic element. Then for any $\mu \in P_+$, $X_{\mu}(\dot x)$ has only finitely many irreducible components. 
\end{cor}

\begin{rmk}
The split case was first proved by Viehmann \cite{Vi} in a different way. 
\end{rmk}

Proof. Let $\pi: \Fl \to Gr$ be the projection. Then $\pi \i X_\mu(\dot x)=\sqcup_{\tw \in W t^\mu W} X_{\tw}(\dot x)$ is a finite union. Since each $X_{\tw}(\dot x)$ has only finite many irreducible components, $\pi \i X_{\mu}(\dot x)$ has only finitely many irreducible components and so is $X_{\mu}(\dot x)$. \qed

\subsection{} Since for superbasic element $b$, $X_{\tw}(b)$ contains only finitely many irreducible components, the number of rational points is finite. Now we show that for split groups, there is a simple formula for this number in terms of class polynomials. 

\begin{prop}
Assume that $G=PGL_n$ is split over $F=\FF_q((\e))$ and $x$ is a superbasic element in $\tW$. Then for any $\tw \in \tW$, $$\sharp X_{\tw}(\dot x)(\FF_q)=n q^{\frac{\ell(\tw)}{2}} f_{\tw, \co} \mid_{v=\sqrt{q}},$$ where $\co$ is the conjugacy class of $\tW$ that contains $x$. 
\end{prop}

Proof. We argue by induction on $\ell(\tw)$. 

Let $\co'$ be the $\d$-conjugacy class that contains $\tw$. 

If $\tw \in \co'_{\min}$, then $\co'=\co$ and $\tw$ is conjugate to $x$. Therefore $\sharp X_{\tw}(\dot x)(\FF_q)=\sharp X_{x}(\dot x)(\FF_q)$. By the equalities (a) and (b) in the proof of Theorem \ref{cct}, $X_{x}(\dot x)(\FF_q) \cong \bJ_{\dot x}/(\bJ_{\dot x} \cap I) \cong Z_{\tW}(x)$.  By the third paragraph of $\S$\ref{5.2}, $Z_{\tW}(x)=\Om$. Hence $\sharp X_{\tw}(\dot x)=n$. 

If $\tw \notin \co'_{\min}$, by Theorem \ref{min} there exists $\tw_1 \approx \tw$ and $i \in \tbS$ such that $\ell(s_i \tw_1 s_i)<\ell(\tw_1)=\ell(\tw)$. By the proof of \cite[Theorem 1.6]{DL}, $$\sharp X_{\tw}(\dot x)(\FF_q)=\sharp X_{\tw_1}(\dot x)(\FF_q)=(q-1) \sharp X_{s_i \tw_1}(\dot x)(\FF_q)+q \sharp X_{s_i \tw_1 s_i}(\dot x)(\FF_q).$$ Hence by inductive hypothesis, 
\begin{align*} \sharp X_{\tw}(\dot x)(\FF_q) &=(q-1) \sharp X_{s_i \tw_1}(\dot x)(\FF_q)+q \sharp X_{s_i \tw_1 s_i}(\dot x)(\FF_q) \\ &=(q-1) n q^{\frac{\ell(\tw)-1}{2}} f_{s_i \tw_1, \co} \mid_{v=\sqrt{q}}+q n q^{\frac{\ell(\tw)-2}{2}} f_{s_i \tw_1 s_i, \co} \mid_{v=\sqrt{q}} \\ &=n q^{\frac{\ell(\tw)}{2}} \bigl( (\sqrt{q}-\sqrt{q} \i) f_{s_i \tw_1, \co} \mid_{v=\sqrt{q}}+f_{s_i \tw_1 s_i, \co} \mid_{v=\sqrt{q}} \bigr) \\ &=n q^{\frac{\ell(\tw)}{2}} f_{\tw, \co} \mid_{v=\sqrt{q}}.
\end{align*} \qed

\section{Reduction method: via partial conjugation}

\subsection{} In this section, we investigate reduction method via ``partial conjugation'' introduced in \cite{He072} and use it to compare the dimension of $X_{\tw}(b)$ for various $\tw$ in the same $W \times W$-coset of $\tW$. 

Notice that any $W \times W$-coset of $\tW$ contains a unique maximal element and this element is of the form $w_0 t^\mu$ for some $\mu \in P_+$. An element in this double coset is of the form $x t^\mu y$ for $x \in W$ and $y \in {}^{I(\mu)} W$. Here $I(\mu)=\{i \in \bS; \<\mu, \a_i\>=0\}$. 

The main result of this section is 

\begin{thm}\label{w0}
Let $\mu \in P_+$, $x \in W$ and $y \in {}^{I(\mu)} W$. Then for any $b \in G(L)$, $$\dim X_{x t^{\mu} y}(b) \le \dim X_{w_0 t^\mu}(b)-\ell(w_0)+\ell(x).$$ In particular, if $X_{x t^\mu y}(b) \neq \emptyset$, then $X_{w_0 t^\mu}(b) \neq \emptyset$. 
\end{thm}

In the rest of this section, we'll prove this Theorem. The proof is rather technical and relies heavily on a detailed analysis of partial conjugation. 
 
\subsection{} We consider the ``partial conjugation'' action of $W$ on $\tW$ defined by $w \cdot_\d w'=w w' \d(w) \i$ for $w \in W$ and $w' \in \tW$. 

For $x \in {}^{\bS} \tW$, set $$I(x)=\max\{J \subset \bS; \Ad(x) \d(J)=J\}.$$ This is well-defined as $\Ad(x) \d(J_1 \cup J_2)=J_1 \cup J_2$ if $\Ad(x) \d(J_i)=J_i$ for $i=1, 2$. 

By \cite[Corollary 2.6]{He072}, we have that $$\tW=\sqcup_{x \in {}^{\bS} \tW} W \cdot_\d (W_{I(x)} x)=\sqcup_{x \in {}^{\bS} \tW} W \cdot_\d (x W_{\d(I(x))}).$$ Moreover, we have the following result (see \cite[Proposition 3.4]{He072}). 

\begin{thm}\label{partial}
For any $\tw \in \tW$, there exists $x \in {}^{\bS} \tW$ and $u \in W_{I(x)}$ such that $$\tw \xrightarrow{i_1}_\d \cdots \xrightarrow{i_k}_\d u x,$$ here $i_1, \cdots, i_k \in \bS$. 
\end{thm}

\

Now we apply this Theorem to the affine Deligne-Lusztig varieties. 

\begin{lem}\label{2}
Let $x, y \in W$ and $\mu \in P_+$ with $y \in {}^{I(\mu)} W$. Then there exists $w' \in {}^{I(\mu)} W$, such that $$\dim X_{x t^\mu y}(b) \le \dim X_{t^\mu w'} (b)+\ell(x).$$
\end{lem}

Proof. We prove the Lemma by induction on $\ell(x)$. Suppose the statement is true for all $x'$ with $\ell(x')<\ell(x)$ but fails for $x$. Let $x=s_{i_1} \cdots s_{i_k}$ be a reduced expression of $x$. There are four different cases

(1) $y s_{\d(i_1)}<y$. In this case,  $y s_{\d(i_1)} \in {}^{I(\mu)} W$.

(2) $y s_{\d(i_1)}>y$ and  $y s_{\d(i_1)} \in {}^{I(\mu)} W$.

(3) $y s_{\d(i_1)}=s_{i_{k+1}} y$ for some $i_{k+1} \in I(\mu)$ and $\ell(s_{i_1} x s_{i_{k+1}})=\ell(x)-2$.

(4) $y s_{\d(i_1)}=s_{i_{k+1}} y$ for some $i_{k+1} \in I(\mu)$ and $\ell(s_{i_1} x s_{i_{k+1}})=\ell(x)$.

If $y s_{\d(i_1)}<y$, then $\ell(x t^\mu y)=\ell(s_{i_1} x t^\mu y s_{\d(i_1)})$. By Proposition \ref{DLReduction}, $\dim X_{x t^\mu y} (b)=\dim X_{s_{i_1} x t^\mu y s_{\d(i_1)}}(b)$. By induction hypothesis for $s_{i_1} x$, there exists $w' \in {}^{I(\mu)} W$ such that $$\dim X_{x t^\mu y} (b) \le \dim X_{t^\mu w'}(b)+\ell(s_{i_1} x)=\dim X_{t^\mu w'}(b)+\ell(x)-1.$$ That contradicts our assumption on $x$.

If $y s_{\d(i_1)}>y$ and  $y s_{\d(i_1)} \in {}^{I(\mu)} W$, then $\ell(x t^\mu y)=\ell(s_{i_1} x t^\mu y s_{\d(i_1)})+2$. By Proposition \ref{DLReduction}, $$\dim X_{x t^\mu y}(b)=\max\{\dim X_{s_{i_1} x t^\mu y s_{\d(i_1)}}(b), \dim X_{s_{i_1} x t^\mu y}(b)\}+1.$$ By induction hypothesis for $s_{i_1} x$, there exists $w' \in {}^{I(\mu)} W$ such that $\dim X_{x t^\mu y} (b) \le \dim X_{t^\mu w'}(b)+\ell(s_{i_1} x)+1=\dim X_{t^\mu w'}(b)+\ell(x)$. That contradicts our assumption on $x$.

Therefore, $y s_{\d(i_1)}=s_{i_{k+1}} y$ for some $i_{k+1} \in I(\mu)$. So $s_{i_1} x t^\mu y s_{\d(i_1)}=s_{i_1} x s_{i_{k+1}}$. If $\ell(s_{i_1} x s_{i_{k+1}})=\ell(x)-2$, then $\ell(x t^\mu y)=\ell(s_{i_1} x s_{i_{k+1}} t^\mu y)+2$.  Hence $$\dim X_{x t^\mu y}(b)=\max\{\dim X_{s_{i_1} x s_{i_{k+1}} t^\mu y}(b), \dim X_{s_{i_1} x t^\mu y}(b)\}+1.$$ By induction hypothesis for $s_{i_1} x$, there exists $w' \in {}^{I(\mu)} W$ such that $\dim X_{x t^\mu y} (b) \le \dim X_{t^\mu w'}(b)+\ell(x)-1$. That contradicts our assumption on $x$.

Hence only case (4) can happen. Now apply the same argument to $s_{i_1} x s_{i_{k+1}}=s_{i_2} \cdots s_{i_{k+1}}$ instead of $x=s_{i_1} \cdots s_{i_k}$, we have that $y s_{\d(i_2)}=s_{i_{k+2}} y$ for some $i_{k+2} \in I(\mu)$ and $\ell(s_{i_2} \cdots s_{i_{k+1}})=\ell(s_{i_3} \cdots s_{i_{k+2}})$. Repeat the same procedure, one may define inductively $i_n \in I(\mu)$ for all $n>k$ by $s_{i_n}=y s_{\d(i_{n-k})} y \i$ and prove that $\ell(s_{i_{n-k+1}} \cdots s_{i_n})=k$. In particular, $i_1, \cdots, i_k \in I(t^\mu y)$. By Proposition \ref{fred}, $\dim X_{x t^\mu y}(b)=\dim X_{t^\mu y}(b)+\ell(x)$. That also contradicts our assumption on $x$. \qed

\subsection{}\label{bed} Now we recall B\'edard's description of ${}^J W$ \cite{B1}. 

Let $\ct(J)$ be the set of all sequences $(J_n, w_n)_{n \ge 1}$, where $J_n \subset \bS$ and $w_n \in W$ such that 

(1) $J_1=J$;

(2) $J_n=J \cap \Ad(w_{n-1}) \d(J_{n-1})$ for $n>1$;

(3) $w_n \in {}^J W^{\d(J_n)}$ and $w_{n+1} \in w_n W_{\d(J_n)}$ for all $n$. 

Then $J_m=J_{m+1}=\cdots$, $w_m=w_{m+1}=\cdots$ and $\Ad(w) \d(J_m)=J_m$ for $m \gg 0$. Moreover $(J_n, w_n)_{n \ge 1} \mapsto w_m$ for $m \gg 0$ is a well-defined bijection between $\ct(J)$ and ${}^J W$. 

\begin{lem}\label{3}
Let $\mu \in P_+$ and $w \in {}^{I(\mu)} W$, then $$\dim X_{t^\mu w}(b) \le \dim X_{w_0 t^\mu}(b)-\ell(w_0).$$
\end{lem}

Proof. Let $J=I(\mu)$ and $(J_n, w_n)_{\ge 1}$ be the element in $\ct(J)$ that corresponds to $w \in {}^J W$. Then there exists $m \in \NN$ such that $J_m=J_{m+1}=\cdots$ and $w_m=w_{m+1}=\cdots$. We prove by descending induction on $n \le m$ that 
\[\tag{a} \dim X_{t^\mu w}(b) \le \dim X_{w_0^{J_n} t^\mu w_n} (b)-\ell(w_0^{J_n}).\]

Notice that $\Ad(w) \d(J_m)=J_m$. Thus $J_m \subset I(t^\mu w)$. By Proposition \ref{fred}, $$\dim X_{w_0^{J_m} t^\mu w} (b)=\dim X_{t^\mu w}(b)+\ell(w_0^{J_m}).$$ Thus (a) holds for $n=m$. Now assume that $n>1$ and (a) holds for $n$, we'll prove that (a) also holds for $n-1$. 

Set $u=\d \i(w_n \i w_{n-1})$. Then $w_{n-1}=w_n \d(u)$. By definition of $\ct(J)$, $u \in W_{J_{n-1}}$, $u \i \in W^{J_n}$ and $\ell(w_n)=\ell(w_{n-1})+\ell(u)$. Now $\ell(u \i w_0^{J_n} t^\mu w_{n-1})=\ell(w_0^{J_n} t^\mu w_n)+2 \ell(u)$. By Proposition \ref{DLReduction}, \[\tag{b}\dim X_{w_0^{J_n} t^\mu w_n}(b) \le \dim X_{u \i w_0^{J_n} t^\mu w_{n-1}}(b)-\ell(u).\]

We prove that 

(c) For any $u' \in W^{J_n}$ and $i \in J_{n-1}$ with $s_i u'<u'$, $\dim X_{u' w_0^{J_n} t^\mu w_{n-1}}(b) \ge \dim X_{s_i u' w_0^{J_n} t^\mu w_{n-1}}(b)+1$.

Note that $w_{n-1} \in {}^{J_{n-1}} W^{\d(J_{n-1})}$. Thus $\ell(w_{n-1} s_{\d(i)})=\ell(w_{n-1})+1$. There are two possibilities.

(i) $w_{n-1} s_{\d(i)} \in {}^{J_{n-1}} W$. 

Then \begin{align*} \ell(s_i u' w_0^{J_n} t^\mu w_{n-1} s_{\d(i)}) &=\ell(s_i u' w_0^{J_n})+\ell(t^\mu)-\ell(w_{n-1} s_{\d(i)}) \\ &=\ell(u' w_0^{J_n})-1+\ell(t^\mu)-\ell(w_{n-1})-1 \\ &=\ell(u' w_0^{J_n} t^\mu w_{n-1})-2. \end{align*}

(ii) $w_{n-1} s_{\d(i)} \notin {}^{J_{n-1}} W$. 

Then $w_{n-1} s_{\d(i)} w_{n-1} \i$ is a simple reflection in $W_{J_{n-1}}$. By $\S$\ref{bed} (2), it is a simple reflection in $W_{J_n}$. Hence \begin{align*} \ell(s_i u' w_0^{J_n} t^\mu w_{n-1} s_{\d(i)}) &=\ell(s_i u' w_0^{J_n} w_{n-1} s_{\d(i)} w_{n-1} \i t^\mu w_{n-1}) \\ &=\ell(s_i u' w_0^{J_n} w_{n-1} s_{\d(i)} w_{n-1} \i)+\ell(t^\mu w_{n-1}) \\ &=\ell(s_i u')+\ell(w_0^{J_n} w_{n-1} s_{\d(i)} w_{n-1} \i)+\ell(t^\mu w_{n-1}) \\ &=\ell(u')+\ell(w_0^{J_n})-2+\ell(t^\mu w_{n-1}) \\ &=\ell(u' w_0^{J_n} t^\mu w_{n-1})-2. \end{align*}

In either case, $\ell(s_i u' w_0^{J_n} t^\mu w_{n-1} s_{\d(i)})=\ell(u' w_0^{J_n} t^\mu w_{n-1})-2$. (c) follows from Proposition \ref{DLReduction}.

Now $u \i \in W_{J_{n-1}} \cap W^{J_n}$. Then $w_0^{J_{n-1}} w_0^{J_n}=v u \i$ for some $v \in W_{J_{n-1}}$ with $\ell(v u \i)=\ell(v)+\ell(u)$. Let $v=s_{i_1} \cdots s_{i_k}$ be a reduced expression, where $k=\ell(w_0^{J_{n-1}} w_0^{J_n})-\ell(u)$. Then $s_{i_j} s_{i_{j+1}} \cdots s_{i_k} u \i \in W^{J_n}$ for all $j$. By (c), 
\begin{align*} \dim X_{u \i w_0^{J_n} t^\mu w_{n-1}}(b) & \le \dim X_{s_{i_k} u \i w_0^{J_n} t^\mu w_{n-1}}(b)-1 \\ & \le \dim X_{s_{i_{k-1}} s_{i_k} u \i w_0^{J_n} t^\mu w_{n-1}}(b)-2 \\ & \le \cdots \\ & \le \dim X_{v u \i w_0^{J_n} t^\mu w_{n-1}}(b)-k \\ &=\dim X_{w_0^{J_{n-1}} t^\mu w_{n-1}}(b)-\ell(w_0^{J_{n-1}} w_0^{J_n})+\ell(u).   \end{align*}

By (b) and induction hypothesis, \begin{align*} \dim X_{t^\mu w}(b) & \le \dim X_{w_0^{J_n} t^\mu w_n} (b)-\ell(w_0^{J_n}) \\ & \le \dim X_{u \i w_0^{J_n} t^\mu w_{n-1}}(b)-\ell(w_0^{J_n})-\ell(u) \\ & \le \dim X_{w_0^{J_{n-1}} t^\mu w_{n-1}}(b)-\ell(w_0^{J_{n-1}} w_0^{J_n})-\ell(w_0^{J_n}) \\ &=\dim X_{w_0^{J_{n-1}} t^\mu w_{n-1}}(b)-\ell(w_0^{J_{n-1}}).
\end{align*}

(a) is proved. 

In particular, $\dim X_{t^\mu w}(b) \le \dim X_{w_0^J t^\mu w_1}(b)-\ell(w_0^J)$. Here $w_1 \in {}^J W^{\d(J)}$. By the same argument as we did for (b), we have that $$\dim X_{w_0^J t^\mu w_1}(b) \le \dim X_{\d \i(w_1) w_0^J t^\mu}(b)-\ell(w_1).$$ 

Similar to the proof of (c), we have that 

(d) Let $u' \in W^J$ and $i \in \bS$ with $s_i u'<u'$. Then $\dim X_{u' w_0^J t^\mu}(b) \ge \dim X_{s_i u' w_0^J t^\mu}(b)+1$. 


Let $w_0 w_0^J \d \i(w_1) \i=s_{j_1} \cdots s_{j_l}$ be a reduced expression, where $l=\ell(w_0)-\ell(w_0^J)-\ell(w_1)$. Then by (d), we have that \begin{align*} \dim X_{w_1 w_0^J t^\mu}(b) & \le \dim X_{s_{i_l} w_1 w_0^J t^\mu}(b)-1 \\ & \le \dim X_{s_{i_{l-1}} s_{i_l} w_1 w_0^J t^\mu}(b)-2 \\ & \le \cdots \\ & \le \dim X_{w_0 t^\mu}(b)-\ell(w_0)+\ell(w_0^J)+\ell(w_1). \end{align*}

Now \begin{align*} \dim X_{t^\mu w}(b) & \le \dim X_{w_0^J t^\mu w_1}(b)-\ell(w_0^J) \\ & \le \dim X_{\d \i(w_1) w_0^J t^\mu}(b)-\ell(w_0^J)-\ell(w_1) \\ & \le \dim X_{w_0 t^\mu}(b)-\ell(w_0).\end{align*} \qed

\subsection{} Now we prove Theorem \ref{w0}. By Lemma \ref{2}, there exists $w \in {}^J W$ such that $\dim X_{x t^\mu y}(b) \le \dim X_{t^\mu w}(b)+\ell(x)$. By Lemma \ref{3}, $\dim X_{t^\mu w}(b) \le \dim X_{w_0 t^\mu}(b)-\ell(w_0)$. Hence $$\dim X_{x t^\mu y}(b) \le \dim X_{w_0 t^\mu}(b)-\ell(w_0)+\ell(x).$$ 

\section{Virtual dimension}

In this section, we discuss some applications of Theorem \ref{w0}. First we compare the dimension of affine Deligne-Lusztig varieties in the affine Grassmannian and the affine flag. 

\begin{thm}\label{compare}
For $\mu \in P_+$ and $b \in G(L)$, $$\dim X_\mu(b)=\dim X_{w_0 t^\mu}(b)-\ell(w_0).$$
\end{thm}

Proof. Let $\pi: \Fl \to Gr$ be the projection. Then each fiber of $\pi$ is isomorphic to $L^+\cg/I$, which is of dimension $\ell(w_0)$. We have that $\pi \i X_\mu(b)=\sqcup_{\tw \in W t^\mu W} X_{\tw}(b)$. 

By Theorem \ref{w0}, $\dim X_{\tw}(b) \le \dim X_{w_0 t^\mu}(b)$ for all $\tw \in W t^\mu W$. Hence 
\begin{align*}
\dim X_{w_0 t^\mu}(b) &=\dim(\sqcup_{\tw \in W t^\mu W} X_{\tw}(b))=\dim \pi \i X_\mu(b) \\ &=\dim X_\mu(b)+\ell(w_0). 
\end{align*}
\qed

\

By combining Theorem \ref{compare} with Theorem \ref{Mazur}, we have the following ``converse to Mazur's inequality'' in the Iwahori case. 

\begin{cor}
Let $J \subset \bS$ with $\d(J)=J$. Let $\mu \in P_+$ and $b \in M_J(L)$ be a basic element such that $\kappa_{M_J}(b) \in Y^+_J$. Then $X_{w_0 t^\mu}(b) \neq \emptyset$ if and only if $\kappa_{M_J}(b) \preceq_J \mu$. 
\end{cor}

\

We define $\eta_\delta\colon \tW\rightarrow W$ as follows. If $\tw=x t^\mu y$ with $\mu \in P_+$, $x \in W$ and $y \in {}^{I(\mu)} W$, then we set $\eta_\delta(x) = \delta^{-1}(y) x$. Then

\begin{thm} 
For any $b \in G(L)$ and $\tw \in \tW$ with $\kappa(\tw)=\kappa(b)$, we have that $$\dim X_{\tw}(b) \le \dim X_\mu(b)+\frac{\ell(\tw)+\ell(\eta_\d(\tw))}{2}-\<\mu, \rho\>.$$

\end{thm}

Proof. We may assume that $\tw=x t^\mu y$ with $x \in W$ and $\mu \in P_+$ and $y \in {}^{I(\mu)} W$. Let $\d \i(y)=s_{i_1} s_{i_2} \cdots s_{i_k}$ be a reduced expression. Since $y \in {}^{I(\mu)} W$, $\ell(t^\mu y)=\ell(t^\mu)-k$ and for $1 \le l \le k$, $t^\mu s_{i_1} \cdots s_{i_l} \in {}^{\bS} \tW$ and $\ell(t^\mu s_{i_1} \cdots s_{i_l})=\ell(t^\mu)-l$. Hence $\ell(\d \i(s_{i_l} \cdots s_{i_1}) \d \i(y) x t^\mu s_{i_1} \cdots s_{i_l})=\ell(\d \i(s_{i_l} \cdots s_{i_1}) \d \i(y) x)+\ell(t^\mu)-l$. 

Notice that $\ell(\d \i(s_{i_l} \cdots s_{i_1}) \d \i(y) x) \ge \ell(\d \i(s_{i_{l-1}} \cdots s_{i_1}) \d \i(y) x)-1$. Thus $$\ell(\d \i(s_{i_{l-1}} \cdots s_{i_1}) \d \i(y) x t^\mu s_{i_1} \cdots s_{i_{l-1}}) \ge \ell(\d \i(s_{i_l} \cdots s_{i_1}) \d \i(y) x t^\mu s_{i_1} \cdots s_{i_l}).$$

So $$\d \i(y) x t^\mu \xrightarrow{i_1}_\d \d \i(s_{i_1}) \d \i(y) x t^\mu s_{i_1} \xrightarrow{i_2}_\d \cdots \xrightarrow{i_k}_\d x t^\mu y.$$ By Proposition \ref{DLReduction}, $$\dim X_{x t^\mu y}(b) \le \dim X_{\d \i(y) x t^\mu}(b)+\frac{1}{2}(\ell(x t^\mu y)-\ell(\d \i(y) x t^\mu)).$$

By Theorem \ref{w0} and Theorem \ref{compare}, \begin{align*} \dim X_{x t^\mu y}(b) & \le \dim X_{w_0 t^\mu}(b)-\ell(w_0)+\frac{1}{2}(\ell(x t^\mu y)+\ell(\d \i(y) x))-\<\mu, \rho\> \\ &=\dim X_\mu(b)+\frac{1}{2}(\ell(x t^\mu y)+\ell(\d \i(y) x))-\<\mu, \rho\>. \end{align*}
\qed

\subsection{} For $\tw\in \tW$ with $\kappa(x)=\kappa(b)$, we define the \emph{virtual dimension}:
\[
d_{\tw}(b)= \frac 12 \big( \ell(\tw) + \ell(\eta_\delta(\tw)) -\defect(b)  \big)-\<\bar \nu_b, \rho\>.
\]

Here $\defect(b)$ is the defect of $b$. See \cite[1.9.1]{Ko3}. 

It is proved in \cite[Theorem 2.15.1]{GHKR1} and \cite[Theorem 1.1]{Vi} that if $G$ is split over $F$, then \[\tag{a} \dim X_\mu(b)=\<\mu-\bar \nu_b, \rho\>-\frac{1}{2} \defect(b).\] As discussed in the proof of Theorem \ref{Mazur}, $\dim X_\mu(b)$ depends only on the combinatorial data $(\tW, \d, \tw, f(b))$. Hence (a) remains true if $\d=id$. Therefore 

\begin{cor}\label{vir}
If $\d=id$, then for any $b \in G(L)$ and $\tw \in \tW$ with $\k(\tw)=\k(b)$, we have that $\dim X_{\tw}(b) \le d_{\tw}(b)$. 
\end{cor}

\section{Class polynomial: lowest two-sided cell case}

We have established an upper bound for the dimension of affine Deligne-Lusztig varieties in the last section. In this section, we study in details of class polynomials for the lowest two-sided cell case. The main result is Theorem \ref{class4}, which is a key step in establishing a lower bound for the dimension of affine Deligne-Lusztig varieties. 

We assume in this section that $G$ splits over $L$. Hence $\Phi=\Sigma$ is a reduced root system. We denote by $\a_i$ (for $i \in \bS$) the set of simple roots. 

\subsection{} 
Let $\tw \in \tW$, let $\supp(\tw)$ be the set of simple reflections of $\tbS$ that appears in some (or equivalently, any) reduced expression of $\tw$ and $\supp_\d(\tw)=\cup_{i \in \NN} \d^i(\supp(\tw))$. 

By \cite[Lemma 1]{H8}, for any $\tw, \tw' \in \tW$, the set $\{ uu';\ u\le \tw,\ u'\le \tw' \}$ has a unique maximal element, which we denote by $\tw \ast \tw'$. Then $\ast$ is associative. Moreover, $\supp(\tw \ast \tw') = \supp(\tw) \cup \supp(\tw')$ and $\supp_\d(\tw \ast \tw') = \supp_\d(\tw) \cup \supp_\d(\tw')$.

Let $A_+=\ZZ_+[v-v \i]$ and $\tH_+=\sum_{\tw \in \tW} A_+ T_{\tw}$. Notice that $T^2_s=(v-v \i)T_s+1$ for any $s \in \tbS$. We have that

\begin{lem}\label{s+}
For any $i \in \tbS$, $\tH_+ T_{s_i} \subset \tH_+$. 
\end{lem}

Proof. Let $\tw \in \tW$. Then \[\tag{a} T_{\tw} T_{s_i}=\begin{cases} T_{\tw s_i}, & \text{ if } \tw s_i>\tw; \\ T_{\tw} T_{s_i}=(v-v \i) T_{\tw}+T_{\tw s_i}, & \text{ otherwise}. \end{cases}\] In particular, $T_{\tw} T_{s_i} \in \tH_+$ and $A_+ T_{\tw} T_{s_i} \subset \tH_+$. \qed

\begin{lem}\label{ttstar}
For any $x, y \in \tW$, we have that 

(1) $T_x T_y \in T_{x y} +\tH_+$. 

(2) $T_x T_y \in (v-v \i)^{\ell(x)+\ell(y)-\ell(x\ast y)} T_{x \ast y}+\tH_+$. 
\end{lem}

Proof. We argue by induction on $\ell(y)$. If $\ell(y)=0$, then $x \ast y=x y$ and $T_x T_y=T_x T_y$. The statement is obvious. 

Now assume that $\ell(y)>0$ and the statement holds for all $y' \in \tW$ with $\ell(y')<\ell(y)$. Let $i \in \tbS$ with $y s_i<y$. Set $y'=y s_i$. Then $y=y' \ast s_i$. By induction hypothesis on $y'$, $$T_x T_y=T_x T_{y'} T_{s_i} \in (T_{x y'}+\tH_+) T_{s_i} \subset T_{x y'} T_{s_i}+\tH_+$$ and 
\begin{align*} 
T_x T_y &=T_x T_{y'} T_{s_i} \in ((v-v \i)^{\ell(x)+\ell(y')-\ell(x \ast y')} T_{x \ast y'}+\tH_+) T_{s_i} \\ & \subset (v-v \i)^{\ell(x)+\ell(y')-\ell(x \ast y')} T_{x y'} T_{s_i}+\tH_+.
\end{align*}

By (a) of Lemma \ref{s+}, $T_{x y'} T_{s_i} \in T_{x y' s_i}+\tH_+$ and $$T_{x \ast y'} T_{s_i} \in (v-v \i)^{\ell(x \ast y')+1-\ell((x \ast y') \ast s_i)} T_{(x \ast y') \ast s_i}+\tH_+.$$ 

Notice that $(x \ast y') \ast s_i=x \ast (y' \ast s_i)=x \ast y$. Hence $T_x T_y \in T_{x y'} T_{s_i}+\tH_+=T_{x y}+\tH_+$ and \begin{align*} T_x T_y & \in (v-v \i)^{\ell(x)+\ell(y')-\ell(x \ast y')} T_{x \ast y'} T_{s_i}+\tH_+ \\ & \subset (v-v \i)^{\ell(x)+\ell(y)-\ell(x\ast y)} T_{x \ast y}+\tH_+.\end{align*} \qed

\

As a consequence, we have that

\begin{cor}\label{th+}
For any $\tw \in \tW$, $T_{\tw} \tH_+ \subset \tH_+$ and $\tH_+ T_{\tw} \subset \tH_+$. 
\end{cor}

\subsection{} We follow \cite[7.3]{Sp}. Let $(W, \bS)$ be a Coxeter system and $\d: W \to W$ which sends simple reflections to simple reflections. For each $\delta$-orbit in $\bS$ we pick a simple reflection and let $c$ be the product of the corresponding simple reflections (in any order). We call $c$ a {\it $\d$-twisted Coxeter element} of $W$. 

The main result we'll prove in this section is 

\begin{thm}\label{class4}
Assume that $G$ is simple. Let $\tW'$ be the lowest two-sided cells of $\tW$. Let $\tw \in \tW' \cap \t W_a$ for $\t \in \Om$. Let $n$ be the number of $\d$-orbits on $\tbS$. Then there exists a maximal proper $\Ad(\t) \circ \d$-stable subset $J$ of $\tbS$ and a $\Ad(\t) \circ \d$-twisted Coxeter element $c$ of $W_J$ such that $$T_{\tw} \in (v-v \i)^{\ell(\eta_\d(\tw))-n} T_{c \t}+\tH_++[\tH, \tH]_\d.$$
\end{thm}

The proof relies on the following three Propositions. 

\begin{prop}\label{xytoa}
Let $\tw=x t^\mu y$ with $\mu \in P_+$, $x \in W$ and $y \in {}^{I(\mu)} W$. If $\tw \in \tW'$ and $\supp_\d(\d \i(y) x)=\bS$, then there exists $a \in W$ with $\supp_\d(a)=\bS$ and $\g \in P_+$ such that $$T_{\tw} \in (v-v \i)^{\ell(\eta_\d(\tw))-\ell(a)} T_{a t^\g}+\tH_+ +[\tH, \tH]_\d.$$
\end{prop}

Proof. Let $J=\{i \in \bS; s_i y<y\}$ and $\rho^\vee_J \in P_+$ with $$\<\rho^\vee_J, \a_i\>=\begin{cases} 1, \text{ if } j \in J \\ 0, \text{ if } j \notin J\end{cases}.$$ Since $y \in {}^{I(\mu)} W$, $J \cap I(\mu)=\emptyset$. Hence $\mu-\rho_J^\vee \in P_+$. Let $J'=I(\mu-\rho^\vee_J)$. Then $\d \i(y)x=w z$ for some $w \in W^{J'}$ and $z \in W_{J'}$. Define $\g \in Y_+$ and $y' \in W^{I(\g)}$ by $\mu-\rho^\vee_J+w \i \rho^\vee_{\d \i(J)} =y' \g$. 

Set $\tw_1=x z \i t^{\mu-\rho^\vee_J} y'$ and $\tw_2=(y') \i z t^{\rho^\vee_J} y$. Then $\tw=\tw_1 \tw_2$ and 
\begin{align*}
\tw_2 \d(\tw_1) & =(y') \i z t^{\rho^\vee_J} y \d(\tw_1)=(y') \i z t^{\rho^\vee_J} y\d(x z \i) \d(t^{\mu-\rho^\vee_J} y')  \\ &=(y') \i z t^{\rho^\vee_J} \d(w) \d(t^{\mu-\rho^\vee_J} y')=(y') \i z \d(w t^{\mu-\rho^\vee_J+w \i \rho^\vee_{\d \i(J)}} y') \\ &=(y') \i z \d(w t^{y' \g} y')=(y') \i z \d(w y' t^\g).
\end{align*}

By \cite[$\S$3.5]{GH}, $\ell(\tw)=\ell(\tw_1)+\ell(\tw_2)$ and $\ell((y') \i z)+\ell(w y')=\ell(\d \i(y) x)$. Set $a =((y')\i z) * \d(w y')$. Then $\supp_\d(a) = \bS$ since $\bS = \supp_\d(\d \i(y)x) \subseteq \supp_\d(w y') \cup \supp_\d((y')\i z)$. We have that 
\begin{align*}
T_{\tw} &=T_{\tw_1} T_{\tw_2} \in T_{\tw_2} T_{\d(\tw_1)}+ [\tH, \tH]_\d=T_{(y') \i z} T_{t^{\rho^\vee_J} y} T_{\d(\tw_1)}+ [\tH, \tH]_\d \\ & \subset T_{(y') \i z} T_{\d(w y' t^\g)}+\tH_++[\tH, \tH]_\d=T_{(y') \i z} T_{\d(w y')} T_{t^\g}+\tH_++[\tH, \tH]_\d \\ &\subset  (v-v \i)^{\ell(\d \i(y)x)-\ell(a)} T_{a} T_{t^\g}+\tH_+ +[\tH, \tH]_\d \\ &\subset  (v-v \i)^{\ell(\d \i(y)x)-\ell(a)} T_{a t^\g}+\tH_+ +[\tH, \tH]_\d.
\end{align*} \qed

\begin{prop}\label{atoc}
Assume that $G$ is simple. Let $J \subset \bS$ with $\d(J)=J$ and $\tw=x t^\mu y$ such that $x \in W_J$ with $\supp_\d(v)=J$, $y$ is a $\d$-twisted Coxeter element in $W_{\bS-J}$, $\mu \neq 0$ and $t^\mu y \in {}^{\bS} \tW$. Then there exists a $\d$-twisted Coxeter element $c$ of $W$ such that $$T_{\tw} \in (v-v \i)^{\ell(x)+\ell(y)-\ell(c)} T_{t^\mu c}+\tH_++[\tH, \tH]_\d.$$ 
\end{prop}

Proof. We proceed by induction on $|J|$. Suppose that the statement is true for all $J'\subsetneq J$, but not true for $J$. We may also assume that the statement is true for all $x'$ with $\supp_\d(x')=J$ and $\ell(x')<\ell(x)$, but not true for $x$. 

Since $G$ is quasi-simple and $\mu \neq 0$, there exists $i \in J$ such that $t^\mu y s_i \in {}^{\bS} \tW$. Let $J_1=\{i \in J; t^\mu y s_i \notin {}^{\bS} \tW\}$. Then $J_1$ is a proper subset of $J$. For any $i \in J_1$, $t^\mu y s_i=s_j t^\mu y$ for some $j \in \bS$. Since $y \in W_{\bS-J}$, this is possible only if $j=i$ and $j$ commutes with $y$. 

We prove that 

(a) $x \in W_{\d \i(J_1)}$. 

We write $x$ as $u x_1$ for $u \in W_{\d \i(J_1)}$ and $x_1 \in {}^{\d \i(J_1)} W$. Then $T_{\tw} \equiv T_{x_1 t^\mu y} T_{\d(u)}=T_{x_1} T_{\d(u)} T_{t^\mu y} \mod [\tH, \tH]_\d$. Let $x'=x_1 \ast \d(u)$. Notice that $\supp_\d(x)=\supp_\d(u) \cup \supp_\d(x_1)=\supp_\d(x_1) \cup \supp_\d(\d(u))=\supp_\d(x')$. Hence $\supp(x')=J$. By Lemma \ref{ttstar} and Corollary \ref{th+}, $$T_{\tw} \in T_{x_1 t^\mu y} T_{\d(u)} +[\tH, \tH]_\d \subset (v-v \i)^{\ell(x)-\ell(x')} T_{x' t^\mu y}+\tH_+ + [\tH, \tH]_\d.$$ 

If $x \notin W_{\d \i(J_1)}$, then $x_1 \neq 1$ and there exists $i \in \d \i(J_1)$ with $s_i x_1<x_1$. Then $s_i x'<x'$. Moreover, $\ell(y s_{\d(i)})=\ell(y)+1$ and $t^\mu y s_{\d(i)} \in {}^{\bS} \tW$. Hence \begin{align*} \ell(s_i x' t^\mu y s_{\d(i)}) &=\ell(s_i x')+\ell(t^\mu)-\ell(y s_{\d(i)})=\ell(x')+\ell(t^\mu)-\ell(y)-2 \\ &=\ell(x' t^\mu y)-2. \end{align*} 

If $\supp_\d(s_i x')=J$, then \begin{align*} T_{x' t^\mu y} &=T_{s_i} T_{s_i x' t^\mu y}\in T_{s_i x' t^\mu y} T_{s_\d(i)}+[\tH, \tH]_\d \\ & \subset (v-v \i) T_{s_i x' t^\mu y}+\tH_+ + [\tH, \tH]_\d.\end{align*} 

By induction hypothesis on $s_i x'$, we have that $$T_{s_i x' t^\mu y} \in (v-v \i)^{\ell(x')-1+\ell(y)-\ell(c)} T_{t^\mu c}+\tH_+ + [\tH, \tH]_\d$$ for some $\d$-twisted Coxeter element $c$ of $W$. 

Therefore $T_{\tw} \in (v-v \i)^{\ell(x)+\ell(y)-\ell(c)} T_{t^\mu c}+\tH_+ + [\tH, \tH]_\d$. That is a contradiction. 

If $\supp_\d(s_i x') \neq J$, then \begin{align*} T_{x' t^\mu y} &=T_{s_i} T_{s_i x' t^\mu y}\in  T_{s_i x' t^\mu y} T_{s_\d(i)}+ [\tH, \tH]_\d \\ & \subset T_{s_i x' t^\mu y s_{\d(i)}}+\tH_+ + [\tH, \tH]_\d.\end{align*} 

By induction hypothesis on $\supp_\d(s_i x')$, we have that $$T_{s_i x' t^\mu y s_i} \in (v-v \i)^{\ell(x')+\ell(y)-\ell(c)} T_{t^\mu c}+\tH_+ + [\tH, \tH]_\d$$ for some $\d$-twisted Coxeter element $c$ of $W$. 

Therefore $T_{\tw} \in (v-v \i)^{\ell(x)+\ell(y)-\ell(c)} T_{t^\mu c}+\tH_+ + [\tH, \tH]_\d$. That is also a contradiction. 

Now (a) is proved. 

We have that $T_{\tw} \in T_{t^\mu y} T_{\d(x)}+ [\tH, \tH]_\d=T_{\d(x) t^\mu y} + [\tH, \tH]_\d$. By the same argument for $\d(x)$ instead of $x$, we have that $\d(x) \in W_{\d \i(J_1)}$. Repeat the same procedure, $\d^i(x) \in W_{\d \i(J_1)}$ for all $i$. Thus $\supp_\d(x) \subset \d \i(J_1) \subsetneqq J$. That is again a contradiction. \qed

\begin{prop}\label{ctoc}
Assume that $G$ is simple. Let $\t \in \Om$. Then there exists a maximal $\Ad(\t) \circ \d$-stable proper subset $J$ of $\tbS$ and a $\Ad(\t) \circ \d$-twisted Coxeter element $c$ of $W_J$ such that $c \t$ is a minimal length element in its $\d$-conjugacy class of $\tW$ and $t^\mu w \to_\d c \t$ for any $\d$-Coxeter element $w$ of $W$ and $\mu \in P$ with $\kappa(t^\mu)=\kappa(\t)$.
\end{prop}

It is proved in a joint work with Yang \cite[Theorem 1.1]{HY}. A partial result for some classical groups was previously obtained in \cite{He5}. 

\subsection{} Now we prove Theorem \ref{class4}. By Proposition \ref{xytoa}, there exists $a \in W$ with $\supp_\d(a)=\bS$ and $\l \in P_+$ such that \[\tag{a}T_{\tw} \in (v-v \i)^{\ell(\eta_\d(\tw))-\ell(a)} T_{a t^\g}+\tH_+ +[\tH, \tH]_\d.\]

By Proposition \ref{atoc}, there exists a $\d$-twisted Coxeter element $w$ of $W$ such that \[\tag{b} T_{a t^\g} \in (v-v \i)^{\ell(a)-n} T_{t^\g w}+\tH_++[\tH, \tH]_\d.\]

Now by Proposition \ref{ctoc}, there exists maximal $\Ad(\t) \circ \d$-stable proper subset $J$ of $\tbS$ and a $\Ad(\t) \circ \d$-twisted Coxeter element $c$ of $W_J$ such that $t^\g a \to_\d c \t$. Thus by $\S$\ref{class}, \[\tag{c}T_{t^\g w} \in T_{c \t}+\tH_++[\tH, \tH]_\d.\]

The theorem then follows from (a), (b) and (c).  

\section{Conjecture of GHKR}

Now we give a lower bound of $\dim X_{\tw}(b)$. 

\begin{thm}\label{upper}
Assume that $G$ is simple. Let $\tw \in \tW'$ and $b \in G(L)$ is a basic element with $\supp_\d(\eta_\d(\tw))=\bS$ and $\kappa(\tw)=\kappa(b)$. Then $\dim X_{\tw}(b) \ge d_{\tw}(b)$. 
\end{thm}

\begin{rmk}
It is proved in \cite{GHKR}, \cite{GH} for a split group and \cite[Theorem B]{GHN} for any tamely ramified group that if $\tw \in \tW'$ and $b \in G(L)$ is a basic element with $\kappa(\tw)=\kappa(b)$, then $X_{\tw}(b) \neq \emptyset$ if and only if $\supp_\d(\eta_\d(\tw))=\bS$. 
\end{rmk}

Proof. By Theorem \ref{class4}, there exists a maximal proper $\Ad(\t) \circ \d$-stable subset $J$ of $\tbS$ and a $\Ad(\t) \circ \d$-twisted Coxeter element $c$ of $W_J$ such that $$T_{\tw} \in (v-v \i)^{\ell(\eta_\d(\tw))-n} T_{c \t}+\tH_++[\tH, \tH]_\d.$$

Let $\co$ be the $\d$-conjugacy class of $\tW$ that contains $c \t$. Then $f_{\tw, \co} \in (v-v \i)^{\ell(\eta_\d(\tw))-n}+A_+$. In particular, $\deg f_{\tw, \co} \ge \ell(\eta_\d(\tw))-n$, where $n$ is the number of $\d$-orbits on $\bS$. 

By Theorem \ref{class3}, \begin{align*} \dim X_{\tw}(b)  & \ge \frac{1}{2}(\ell(\tw)+\ell(c)+\deg f_{\tw, \co}) \\ & \ge \frac{1}{2}(\ell(\tw)+\ell(\eta_\d(\tw))+\ell(c)-n). \end{align*}

By the definition of defect, $\defect(b)=\defect(\dot \t)=n-\ell(c)$. Therefore $\dim X_{\tw}(b) \ge d_{\tw}(b)$. \qed 

\

By combining Theorem \ref{upper} and Corollary \ref{vir}, we have that 

\begin{cor}
If $G$ is simple and $\d=id$. Let $\tw \in \tW'$ and $b \in G(L)$ be a basic element with $\supp(\eta(\tw))=\bS$ and $\kappa(\tw)=\kappa(b)$. Then $\dim X_{\tw}(b)=d_{\tw}(b)$. 
\end{cor}

\begin{rmk}
The split case was first conjectured in \cite[Conjecture 1.1.3]{GHKR}. A weaker result for some split classical groups was proved in \cite{GH}. 
\end{rmk}

\section*{Acknowledgment} We thank G. Lusztig, U. G\"ortz and T. Haines for useful comments. We also thank the referee for careful reading and useful suggestions.

\end{document}